\documentclass[12pt]{article}

\usepackage{xy} \xyoption{all}
\usepackage{hyperref}
\usepackage{amsthm}
\usepackage{amsmath}
\usepackage{amssymb}
\usepackage{amsfonts}
\usepackage{ifthen}

\def\Fp{{{\mathbb F}_p}}

\def\QQ{{\mathbb Q}}

\def\ZZ{{\mathbb Z}}

\def\urltilde{\kern-.01em\lower.6ex\hbox{\~{}}}
 \def\p{\mathfrak p}

\def\Norm{N} \def\q{{\mathfrak q}}
\def\McOmega{K_S} \def\McG{\mathcal{G}} \def\McR{\mathcal{R}}
\def\McF{\mathcal{F}}   
\def\g{\mathfrak g}

\def\Mcrho{\rho} \def\Mczeta{\zeta} \def\McO#1{\mathcal{O}_{{#1},S}}
\def\McU#1{\McO{#1}^{\times}} 
\def\maxp{10,\!000} \def\zp{\ZZ_p}
\def\McD{\mathcal{D}} \def\ab{\mathrm{ab}}
\def\cl#1{C_{{#1},S}} \def\clnos#1{C_{#1}} \def\pcl#1{A_{{#1},S}}
\def\Mcpairing#1#2{({#1},{#2})_{S}}
\def\blankpairing{\Mcpairing{\ }{\ }}
\def\eigpair#1#2#3{\langle{#1},{#2}\rangle_{#3}}
\def\blankeigpair#1{\eigpair{\ }{\ }{#1}}

\def\McOKS{\mathcal{O}_S} \def\McUKS{\McOKS^{\times}}
\def\Fil{\mathrm{Fil}}
\def\McI#1{I_{{#1},S}} \def\McP#1{P_{{#1},S}} \def\McH#1{H_{{#1},S}}
\def\McA#1{A_{{#1},S}}
\def\h{\mathfrak{h}}

\def\into{\hookrightarrow}
\def\iso{\simeq}
\def\to{\rightarrow}
\def\tensor{\otimes}
\def\divides{\mid}

  \def\Aut{\operatorname{Aut}}
  \def\Gal{\operatorname{Gal}}
  \def\gr{\operatorname{gr}}
  \def\Hom{\operatorname{Hom}}
  \def\inv{\operatorname{inv}}
  \def\ord{\operatorname{ord}}
  \def\res{\operatorname{res}}
  \def\cor{\operatorname{cor}}
  \def\Fil{\operatorname{Fil}}

    \ifthenelse{\value{secnumdepth}=0}{\newtheorem{theorem}{Theorem}}
    {\newtheorem{theorem}{Theorem}[section]}
    \newtheorem{lemma}[theorem]{Lemma}
    \newtheorem{proposition}[theorem]{Proposition}
    \newtheorem{corollary}[theorem]{Corollary}
    \newtheorem{conjecture}[theorem]{Conjecture}

\ifthenelse{\value{secnumdepth}=0}{}{\setcounter{secnumdepth}{2}}

\newcommand{\etalchar}[1]{$^{#1}$}
\providecommand{\bysame}{\leavevmode\hbox
to3em{\hrulefill}\thinspace}
\providecommand{\MR}{\relax\ifhmode\unskip\space\fi MR }
 \providecommand{\href}[2]{#2}

\title{A cup product in the Galois cohomology\\ of number fields}

\author{William G. McCallum and Romyar T. Sharifi}

\begin{document}
\maketitle

\section{Introduction} \label{intro}
This paper is devoted to the consideration of a certain cup
product in the Galois cohomology of an algebraic number field with
restricted ramification.  Let $n$ be a positive integer, let $K$
be a number field containing the group $\mu_{n}$ of $n$th roots of
unity, and let $S$ be a finite set of primes including those above
$n$ and all real archimedean places.  Let $G_{K,S}$ denote the
Galois group of the maximal extension of $K$ unramified outside
$S$ (inside a fixed algebraic closure of $K$).  We consider the
cup product
\begin{equation} \label{cupprod}
  H^1(G_{K,S},\mu_{n}) \otimes H^1(G_{K,S},\mu_{n}) \to
  H^2(G_{K,S},\mu_{n}^{\otimes 2}).
\end{equation}
When the localization map
\begin{equation}\label{localizationmap}
H^{2}(G_{K,S}, \mu_{n}^{\otimes 2}) \to \bigoplus_{v \in S}
H^{2}(G_{v}, \mu_{n}^{\otimes 2})
\end{equation}
for the local absolute Galois groups $G_{v}$ is injective, the cup
product is a direct sum over primes in $S$ of the corresponding
local cup products, each of which may be expressed as the $n$th
norm residue symbol on the completion of $K$ at that prime.
However, we are interested in a situation that is inherently
non-local: that is, in which the cup product of two elements lies
in the kernel of the localization map.  This kernel is isomorphic
to $\cl{K}/n\cl{K}$, where $\cl{K}$ denotes the $S$-class group of
$K$. In Section~\ref{cupproduct}, we develop a formula for the cup
product in this case in terms of ideals in a Kummer extension of
$K$ (Theorem~\ref{PairingFormula}).

Using Kummer theory, we can identify $H^1(G_{K,S},\mu_{n})$ with a
subgroup of $K^{\times}/K^{\times n}$ containing the image of the
units of the $S$-integers $\McO{K}$, and hence we have an induced
pairing
\[  \McU{K} \times \McU{K} \to H^2(G_{K,S},\mu_{n}^{\otimes 2}). \]
Like the norm residue symbol, this pairing has the property that $a$
and $b$ pair trivially if $a + b = 1$.  Thus, we obtain a relationship
between the cup product and the $K$-theory of $\McO{K}$, which is
described in Section~\ref{K-theory} and discussed further in
Section~\ref{cyclpunits} (see Conjecture~\ref{k2conjecture}).

When $n = p$, a prime number, the cup product yields information
on the form of relations in the maximal pro-$p$ quotient $\McG =
G_{K,S}^{(p)}$ of $G_{K,S}$. It is well-known that $\McG$ has a
presentation
\begin{equation} \label{presentation}
0 \to \McR \to \McF \to \McG \to 0,
\end{equation}
where $\McF$ is a free pro-$p$ group on a finite set $X$ of
generators and $\McR$ is the smallest closed normal subgroup
containing a finite set $R$ of relations in $\McF$. Choosing $X$
and $R$ of minimal order, we have
$$
|X| = \dim_{\ZZ/p\ZZ} H^{1}(\McG, \ZZ/p\ZZ) \quad\text{and}\quad
|R| = \dim_{\ZZ/p\ZZ}H^{2}(\McG, \ZZ/p\ZZ).
$$
When the localization map \eqref{localizationmap} is injective, the
relations can be understood in terms of relations in decomposition
groups.  On the other hand, it is quite difficult to say anything
about relations corresponding to the kernel of the localization map.

Let us quickly review the precise relationship between
$H^{i}(\McG,\ZZ/p\ZZ)$, $i = 1, 2$, and the generators and
relations for $\McG$, as detailed in \cite{labute:1967} or
\cite[Section~III.9]{neukirchetal:2000}.  Let $\gr^{\cdot}(\McF)$
denote the sequence of graded quotients associated with the
descending $p$-central series on $\McF$.  The image of $X$ in
$\gr^{1}(\McF)= \McF/\McF^{p}[\McF,\McF]$ forms a basis, which we
also denote by $X$.  With the choice of a linear ordering on $X$,
the set
\[ \{p x, [x, x']\colon x, x' \in X,\ x
< x'\} \] forms a basis for $\gr^{2}(\McF)$, where $p:
\gr^{1}(\McF) \to \gr^{2}(\McF)$ is induced by the $p$th power map
and $[x,x']$ is (the image of) the commutator
$xx'x^{-1}(x')^{-1}$. Now, any relation $\Mcrho \in R$ has zero
image in $\gr^{1}(\McF)$, and its image in
$\gr^{2}(\McF)/p\gr^{1}(\McF)$ is
\begin{equation}
     \sum_{x < x' \in X}
     a^{\Mcrho}_{x,x'}[x,x'], \quad
a^{\Mcrho}_{x,x'} \in
     \ZZ/p\ZZ.\label{rho}
\end{equation}
The quotient $\gr^{1}(\McF)$ and \( H^{1}(\McG, \ZZ/p\ZZ) \iso
H^1(\McF,\ZZ/p\ZZ) \) are dual as vector spaces over $\Fp$,
allowing us to define a basis $X^{*}$ of $H^{1}(\McG, \ZZ/p\ZZ)$
dual to $X$. Similarly, $R$ may be regarded as a basis for the
dual to $H^{2}(\McG, \ZZ/p\ZZ)$ via the transgression isomorphism
$H^1(\McR,\mathbb{Z}/p\mathbb{Z})^{\McG} \iso
H^2(\McG,\mathbb{Z}/p\mathbb{Z})$.  With these identifications,
the numbers $a^{\Mcrho}_{x, x'}$ in \eqref{rho} are given by
\begin{equation}
a^{\Mcrho}_{x, x'} = -\rho(x^{*} \cup x^{\prime *}), \quad x <
x',\label{BasicRelation}
\end{equation}
where the cup refers to the cup product pairing
\begin{equation} \label{zmodpzcup}
H^{1}(\McG, \ZZ/p\ZZ) \times H^{1}(\McG, \ZZ/p\ZZ) \to H^{2}(\McG,
\ZZ/p\ZZ).
\end{equation}
Since $K$ contains $\mu_p$, we have natural isomorphisms
\[  H^i(G_{K,S},\mu_{p}^{\otimes j}) \iso H^i(\McG,\ZZ/p\ZZ) \otimes
\mu_{p}^{\otimes j} \] for any $i$ and $j$, and so the cup product
\eqref{zmodpzcup} is just \eqref{cupprod} after a choice of
isomorphism $\mu_p \iso \ZZ/p\ZZ$.

Our primary focus in
Sections~\ref{IwasawaLfunction}--\ref{idealclass} is a case in
which the localization map \eqref{localizationmap} is zero: $n =
p$, with $p$ an odd prime, $K = \QQ(\mu_{p})$, and $S$ consisting
of the unique prime above $p$.  In this case there is a natural
conjugation action of $\Delta = \Gal(K/\QQ)$ on $\gr^{i}(\McF)$;
let us suppose that generators and relations have been chosen so
that each $x$ and each $\Mcrho$ is an eigenvector for this action.
Define $a_{x,x'}^{\Mcrho}$ by equation \eqref{BasicRelation} for
any $x,x'\in X$.  Suppose further that $p$ satisfies Vandiver's
conjecture.  Then $H^{2}(\McG,\ZZ/p\ZZ)^{-}$, the subspace on
which $\Delta$ acts by an odd character, is trivial, and therefore
the matrix $A^{\Mcrho} = (a_{x,x'}^{\Mcrho})$ decomposes into
blocks
\begin{equation*}
     A^{\Mcrho} =
    \left(
    \begin{array}{cc}
        A^{\Mcrho+} &  0 \\
        0 & A^{\Mcrho-}
     \end{array}\right)
\end{equation*}
where $A^{\Mcrho\pm}$ contains all entries $a_{x,x'}^{\Mcrho}$
corresponding to $x,x' \in \gr^{1}( \McF )^{\pm}$. The matrix
$A^{\Mcrho+}$ is related to the $p$-adic zeta-function (see
Sect.~\ref{IwasawaLfunction}) and can be shown to be nonzero when
Vandiver's conjecture holds and the $\lambda$-invariant of the
$p$-part of the class group in the cyclotomic tower is equal to
its index of irregularity (Proposition~\ref{zetapairing}), in
particular for $p < 12,\!000,\!000$
\cite{buhler-crandall-ernvall-metsankyla-shokrollahi-2001}.

The question of when $A^{\Mcrho-}$ is non-zero is more mysterious.
In Section~\ref{cyclpunits}, we present a method for imposing
linear conditions on $A^{\Mcrho-}$. The relations $\rho \in R$
correspond to nontrivial eigenspaces of the $p$-part of the class
group of $K$.
For all such eigenspaces for $p$ with $p < \maxp$, the method
specifies $A^{\Mcrho-}$ up to a scalar multiple
(Theorem~\ref{uniquepairing}). Since the method works by imposing
linear conditions, it is not capable of showing that
\begin{equation}
     A^{\Mcrho-}\neq 0.\label{Akneq0}
\end{equation}

In Sections~\ref{idealclass} and~\ref{AnotherFormula}, we consider
conditions for
the nontriviality in \eqref{Akneq0}.
In Section~\ref{idealclass}, we describe the
relation between
this nontriviality
and the structure of class groups of Kummer
extensions of $K$. In Section~\ref{AnotherFormula}, we determine a
formula for a certain projection of the cup product when it can be
expressed as the corestriction of a cup product in an unramified
Kummer extension $L$ of $K$ (Theorem~\ref{HCFpairformula}).  We
then describe a
computer
calculation
that uses this formula to verify the nontriviality of
$A^{\Mcrho-}$ for $p = 37$ (Theorem~\ref{NontrivialPairing}).

In Section~\ref{GaloisGroup}, we describe, in detail, the
relations in $\McG$ that result from the nontriviality in
\eqref{Akneq0} (Theorem~\ref{grouprelation}) and, in particular,
the explicit relation for $p = 37$.  In Section~\ref{fundamental},
we use this description to exhibit relations in a graded $\zp$-Lie
algebra $\mathfrak{g}$ associated with the action of
the absolute Galois group
on the
pro-$p$ fundamental group of ${\mathbb P}^1 - \{0,1,\infty\}$
(Theorem~\ref{lierelation}).  For $p = 691$, we explain how a
conjecture of Ihara on the relationship between the structure of
$\mathfrak{g}$ and a certain Lie algebra of derivations implies
\eqref{Akneq0} and, conversely, our calculations confirm Ihara's
conjecture in this case if \eqref{Akneq0} is satisfied
(Theorem~\ref{p691m12}).

In Section~\ref{greenberg}, we consider an Iwasawa-theoretic
consequence of the nontriviality in \eqref{Akneq0}, namely,
Greenberg's pseudo-null conjecture in the case that the $p$-part of
the class group of $K$ has order $p$ (Theorem~\ref{greenbergconj}).

\smallskip

{\bf Acknowledgments.} The authors would like to thank Ralph
Greenberg, Yasutaka Ihara, and Barry Mazur for their support,
encouragement, and useful comments, Robert Bond and Michael Reid
for some corrections, John Cremona for his invaluable assistance
with C++ (used to prove Theorem~\ref{uniquepairing} up to $6500$),
Claus Fieker for determining the structure of the maximal order of
a degree $37$ number field as part of the computation for $p=37$
of Section~\ref{AnotherFormula}, William Stein for his substantial
contribution to this computation, and Kay Wingberg for his
suggestion of the approach in Section~\ref{AnotherFormula}.  The
first author was supported by NSF grant DMS-9624219 and by an AMS
Bicentennial Fellowship and the second author by NSF VIGRE grant
9977116 and by an NSF Postdoctoral Research Fellowship.

\section{A formula for the cup product} \label{cupproduct}

We consider the general setting for the cup product
\eqref{cupprod} of the introduction. First, we describe the groups
$H^{i}(G_{K,S}, \mu_{n})$, $i = 1,2$. Let $\McOmega$ denote the
maximal extension of $K$ unramified outside $S$.  From the
inflation-restriction exact sequence, $H^{1}(G_{K,S}, \mu_{n})$ is
the kernel of the restriction map from $H^{1}(K, \mu_{n})$ to
$H^{1}(\McOmega, \mu_{n})$ and hence may be identified using
Kummer theory with the kernel of $K^{\times}/K^{\times n} \to
\McOmega^{\times}/\McOmega^{\times
  n}$. Thus
$$
  H^{1}(G_{K,S}, \mu_{n}) \iso D_K/K^{\times n},
$$
where
$$
  D_K = \McOmega^{\times n} \cap K^{\times} =
  \{x \in K^{\times}: \mbox{$n \divides
  \ord_{\mathfrak{q}}(x)$ for all $\mathfrak{q} \notin S$}\}.
$$

The description of $H^{2}(G_{K,S}, \mu_{n})$ involves the
$S$-ideal class group $\cl{K}$ of $K$.  For an extension $F/K$,
denote by $\McO{F}$ the ring of $S$-integers in $F$. For brevity,
we set $\McOKS = \McO{\McOmega}$. Since $\McOmega$ contains the
Hilbert class field of $K$, any nonzero ideal $\mathfrak{a}$ of
$\McO{K}$ becomes principal in $\McOKS$---say $\mathfrak{a} =
(\alpha)$ for $\alpha \in \McOmega^{\times}$. Furthermore, since
$\mathfrak{a}$ is fixed by $G_{K,S}$, we have
$\alpha^{\sigma}/\alpha \in \McUKS$ for $\sigma \in G_{K,S}$, and
thus we can associate with $\mathfrak{a}$ an element of
$H^{1}(G_{K,S}, \McUKS)$.  Using Hilbert's Theorem 90 and the
exact sequence
$$
1 \to \McUKS \to \McOmega^{\times} \to P_S \to 1,
$$
where $P_S$ is the group of principal ideals of $\McOKS$, it is
easy to see that this induces an isomorphism
\begin{equation*}
     \cl{K}\iso H^{1}(G_{K,S},\McUKS).
\end{equation*}
Also, if $I_S$ is the group of fractional ideals of $\McOKS$,
then, since $\McOmega$ contains the Hilbert class field of any of
its subfields, $I_S/P_S$ is trivial.  Hence, since $I_S$ is a
direct sum of induced modules, we have
\[
  H^{i}(G_{K,S}, P_S) \iso H^{i}(G_{K,S}, I_S) = 0, \quad i \geq 1.
\]
It follows that $H^{2}(G_{K,S},\McUKS)$ is isomorphic to
$H^{2}(G_{K,S},\McOmega^{\times})$, the $p$-part of which, for
those $p$ dividing $n$, may be identified with the $p$-part of the
subgroup of the Brauer group of $K$ consisting of elements with
zero invariant at all valuations $v \notin S$
\cite[Proposition~8.3.10]{neukirchetal:2000}.  As $S$ contains all
primes dividing $n$, the sequence
\begin{equation} \label{kummer}
1 \to \mu_{n} \to \McUKS \xrightarrow{n}{} \McUKS \to 1
\end{equation}
is exact.  Taking its cohomology and tensoring with $\mu_{n}$, we
obtain an exact sequence
\begin{equation}
  1 \to \cl{K}/n \cl{K} \otimes \mu_{n} \to H^{2}(G_{K,S},
  \mu_{n}^{\otimes 2})
  \xrightarrow{\pi} \bigoplus_{v \in S} \mu_{n}
  \xrightarrow{\prod} \mu_{n}
  \to 1, \label{H2isAmodp}
\end{equation}
where $\pi$ is the twist by $\mu_{n}$ of the direct sum of the
invariant maps.  For future reference, we record that if the ideal
$\mathfrak{a}$ represents an element of $\cl{K}/n \cl{K}$, the
corresponding element of $H^{2}(G_{K,S}, \mu_{n})$ is the
coboundary $H^{1}(G_{K,S}, \McUKS) \to H^{2}(G_{K,S}, \mu_{n})$ of
the cocycle $\alpha^{\sigma}/\alpha$, where $\mathfrak{a}\McOKS =
(\alpha)$.

Now we consider the pairing
$$
  \blankpairing = (\ ,\ )_{n,K,S} \colon D_K \times D_K \to
  H^{2}(G_{K,S}, \mu_{n}^{\otimes 2})
$$
induced by the cup product. Let $a,b \in D_K$.  The image of
$\Mcpairing{a}{b}$ under $\pi$ as in \eqref{H2isAmodp} is given by
the direct sum of the Hilbert pairings $(a,b)_v$ at $v \in S$. We
determine a formula for $\Mcpairing{a}{b}$ with $a,b \in D_K$ in
the case that $\pi\Mcpairing{a}{b}$ is trivial, obtaining an
element of $\cl{K}/n\cl{K} \tensor \mu_n$.

Fix a primitive $n$th root of unity $\zeta$.
Let $\alpha \in K_S^{\times}$ with $\alpha^n = a$, and for $\sigma
\in G_{K,S}$, let $m_{\sigma}$ denote the smallest nonnegative
integer such that $\sigma \alpha = \zeta^{m_\sigma} \alpha$.   We
start with a general cohomological lemma.

\begin{lemma} \label{minus1}
  Consider the homomorphism
  $G_{K,S} \to \ZZ/n\ZZ$ associated with $a$ via Kummer theory and
  the isomorphism $\mu_n \iso \ZZ/n\ZZ$ given by $\zeta \mapsto 1$.
  Let $\epsilon$ denote
  its coboundary in the cohomology sequence
of
  \begin{equation*}
    1 \to \ZZ/n\ZZ \to \ZZ/n^2\ZZ \to \ZZ/n\ZZ \to 1.
  \end{equation*}
  Then
  $$
    \Mcpairing{a}{a} = \frac{n(n-1)}{2}\epsilon \tensor \zeta^{\otimes 2}.
  $$
\end{lemma}

\begin{proof}
  This is antisymmetry of the cup product if $n$ is odd, since in that
  case both sides are zero.  For $n$ even,
  one can check directly that the difference of cocycles
  $$
    (\sigma,\tau) \mapsto m_{\sigma}m_{\tau} -
    (m_{\sigma}+m_{\tau}-m_{\sigma\tau})/2 \pmod{n}
  $$
  is the coboundary of
  $$
    \sigma \mapsto -m_{\sigma}(1+m_{\sigma})/2 \pmod{n}.
  $$
  (See also \cite[Section~III.9]{neukirchetal:2000}.)
\end{proof}

We proceed by considering the exact sequence of $G_{K,S}$-modules
\begin{equation}\label{eq:SLn}
  1 \to \mu_{n} \to SL_{n}(\McOmega) \to PSL_{n}(\McOmega) \to
  1.
\end{equation}
We define a $G_{K,S}$-cocycle with values in $PSL_n(\McOmega)$,
the coboundary of which, when computed in two different ways,
yields the formula. (The choice of this cocycle was motivated by
consideration of the cyclic algebra associated with $a$ and $b$.)

Let
\begin{equation*}
  M_b=
  \left[
  \begin{array}{cccc}
    0 & 1 & \dots & 0\\
    \vdots & \ddots & \ddots & \vdots\\
    0 & \dots & 0 & 1\\
    b & 0 & \dots & 0
  \end{array}
  \right].
\end{equation*}
Note that $M_b^n = b$, and choose $\beta \in \McOmega^{\times}$
such that
$$
  {\beta}^n = \det M_b = (-1)^{n-1}b.
$$
We define a cochain with values in $SL_n(\McOmega)$ by
\begin{equation*}
  {C}^{a,b}_\sigma = (\beta^{-1}M_b)^{m_\sigma}, \quad
  \sigma \in G_{K,S}.
\end{equation*}
The images of the $C_{\sigma}^{a,b}$ in $PSL_n(\McOmega)$ define a
cocycle, as each is fixed by $G_{K,S}$ and has order dividing $n$.
We denote the class of this cocycle by $C^{a,b}$.

\begin{lemma}\label{cupandazumaya}
  Let $\delta$ denote the coboundary in cohomology associated with
  \eqref{eq:SLn}.  Then
  \begin{equation*}
    \delta(C^{a,b}) \otimes \zeta = \Mcpairing{a}{-ab}.
  \end{equation*}
\end{lemma}

\begin{proof}
  Let $n_\sigma$ be the smallest nonnegative integer such that $\sigma
  \beta = \zeta^{n_\sigma} \beta$.  The coboundary of ${C}^{a,b}$ in
  $H^2(G_{K,S}, \mu_n)$ is represented by the 2-cocycle
  \begin{eqnarray*}
    (\sigma, \tau)&\mapsto&(\beta^{-1} M_b)^{m_\sigma}(\beta^{-1}
    M_b)^{m_\tau\sigma}(\beta^{-1} M_b)^{-m_{\sigma\tau}}\\
    && = \beta^{-(\sigma-1)m_\tau}
    (\beta^{-1}M_b)^{{m_\sigma + m_\tau - m_{\sigma\tau}}}\\
    && =  \zeta^{-n_\sigma m_\tau}
    (-1)^{(n-1)(m_\sigma + m_\tau - m_{\sigma\tau})/n},
  \end{eqnarray*}
  which, upon tensoring with $\zeta$ and applying Lemma \ref{minus1},
  yields
  \begin{equation*}
    -\Mcpairing{(-1)^{n-1}b}{a}+ \Mcpairing{a}{a} =
    \Mcpairing{a}{-ab}.
  \end{equation*}
\end{proof}

The following lemma is applicable to the case of interest that
$\pi\Mcpairing{a}{b}$ is trivial.

\begin{lemma}\label{thm:lemma-let-b}
  Let $B$ be the class of a $G_{K,S}$-cocycle $B_{\sigma}$ with
  values in $PSL_n(\McOmega)$ such that
  $\pi(\delta(B) \otimes \zeta) = 0$.
  Then there exists an $A \in GL_n(\McOmega)$
  such that
  $$
    A^\sigma A^{-1} \equiv B_\sigma \pmod{\McOmega^\times}.
  $$
  For any such $A$, there is a fractional ideal
  $\mathfrak{a}$ of $\McO{K}$ such that
  $$
    \det(A)\mathcal{O}_S \equiv \mathfrak{a}\mathcal{O}_S
    \pmod{nI_S}
  $$
  and
  $$\delta(B) = -\mathfrak{a} \pmod{nC_{K,S}}.$$
\end{lemma}

\begin{proof}
  Consider the commutative diagram
  $$
    \xymatrix{
    &1 \ar[d] & 1 \ar[d] & 1 \ar[d] & \\
    1 \ar[r] & \mu_{n} \ar[d]^f \ar[r] & SL_{n}(\McOmega) \ar[r]
    \ar[d] & PSL_{n}(\McOmega) \ar[r]
    \ar[d]^g  & 1\\
    1 \ar[r] & \McOmega^\times \ar[r] \ar[d]& GL_{n}(\McOmega)
    \ar[r] \ar[d]^{\mathrm{det}}&
    PGL_{n}(\McOmega) \ar[r] \ar[d]& 1\\
    1\ar[r] & K_S^{\times n} \ar[d] \ar[r] & \McOmega^\times
    \ar[r]\ar[d] &
    \McOmega^\times/\McOmega^{\times n} \ar[d]\ar[r] & 1\\
    & 1 & 1 & 1 & \\}
  $$
  We use $\delta$ (resp., $\delta'$)
  to denote any coboundary in a long exact sequence of
  cohomology groups arising from a horizontal (resp., vertical)
  short exact sequence in this diagram.
The map $f_* \tensor \mu_n$, with
$$
 f_*: H^2(G_{K,S}, \mu_n) \rightarrow
 H^2(G_{K,S},\McOmega^\times),
$$
can be identified with $\pi$, so our assumption on $B$ is
equivalent to $f_*(\delta(B)) = 0$.  Note that we have
$\delta(g_*(B)) = f_*(\delta(B))$, where $g_*$ denotes the map
$$
  g_* \colon H^1(G_{K,S}, PSL_n(\McOmega)) \to H^1(G_{K,S},
  PGL_n(\McOmega)).
$$
Since $H^1(G_{K,S}, GL_n(\McOmega)) = 0$, this implies $g_*(B) =
0$.  Hence, there exist an element $y \in
(\McOmega^\times/\McOmega^{\times n})^{G_{K,S}}$ such that $B =
\delta'(y)$ and so, by definition of $\delta'$, an $A \in
GL_n(K_S)$ satisfying both the first statement of the lemma and
\begin{equation}
y = \det(A) \pmod{\McOmega^{\times n}}.\label{eq:2}
\end{equation}
(Note that it suffices to check the
remainder of the lemma for this choice of $A$.)

Since the coboundaries $\delta$ and $\delta'$ anticommute,
  $$
    \delta(B) = \delta(\delta'(y)) = -\delta'(\delta(y)).
  $$
Let $\mathfrak{a}$ be a fractional ideal of $\mathcal{O}_{K,S}$
such that
\begin{equation}
\mathfrak{a}\mathcal{O}_S \equiv y\mathcal{O}_S
\pmod{nI_S}.\label{eq:4}
\end{equation}
Choose $\eta \in K_S^{\times}$ such that $\eta\mathcal{O}_S =
\mathfrak{a}\mathcal{O}_S $, and observe that $y = \eta
\pmod{K_S^{\times n}}$. Thus, $\delta(y)$ is the image of the
cocycle associated to $\eta$ under the map
$$
  H^1(G_{K, S},
  \mathcal{O}_S^\times) \rightarrow H^1(G_{K, S}, K_S^{\times n}),
$$
so $\delta'(\delta(y))$ can be computed using the Kummer
sequence~\eqref{kummer}. Hence, by the discussion following
\eqref{H2isAmodp}, $\delta'(\delta(y))$ is the class modulo $n$ of
$\mathfrak{a}$, which, by~\eqref{eq:2} and~\eqref{eq:4}, satisfies
the conditions of the lemma.
\end{proof}

Let $\Norm_{L/K}$ denote the norm map for an extension $L$ of $K$.
We denote the image of a (fractional) ideal $\mathfrak{a}$ in
$\cl{K}/n \cl{K}$ by $[\mathfrak{a}]$.  We remark that
$(n/2)[\mathfrak{a}]$ is always trivial if $n$ is odd.  We now
state our formula for the pairing.

\begin{theorem}\label{PairingFormula}
  Let $a, b \in D_K$ be such that the Hilbert pairing $(a,b)_{v}$ is
  trivial for all valuations $v \in S$.  Choose $\alpha \in
  \McOmega^{\times}$ such that $\alpha^{n}=a$, let $L =
  K(\alpha)$, and set $d = [L:K]$.
  Let $\mathfrak{b}$ be the
  fractional ideal of
  $\McO{K}$ such that $b\McO{K} = \mathfrak{b}^n$.  Write $b=
  \Norm_{L/K}\gamma$ for some $\gamma$, and write
  \begin{equation} \label{yfactorization}
    \gamma\McO{L}= \mathfrak{c}^{1-\sigma}\mathfrak{b}^{n/d},
  \end{equation}
  for some fractional ideal $\mathfrak{c}$ of $\McO{L}$ and $\sigma \in
  \Gal(L/K)$.  Let $\xi \in \mu_n$ be such that $\sigma \alpha =
  \xi \alpha$.  Then
  $$
    \Mcpairing{a}{b} = ([\Norm_{L/K}(\mathfrak{c})]+ \frac{n}{2}
    [\mathfrak{b}])\tensor \xi.
  $$
\end{theorem}

\begin{proof}
  We have that $b$ is a norm from $L_v$ for all $v$, since the
  Hilbert pairing $(a,b)_{v}$ is trivial for $v \in S$ by assumption and
  for $v \notin S$ by definition of $D_K$.
  Since, in a cyclic extension, an element that is a local norm
  everywhere is a global norm, we have $b = \Norm_{L/K} \gamma$ with
  $\gamma \in L^{\times}$.  It is easy to see that $\gamma$ has a
  decomposition as in \eqref{yfactorization}.  Without loss of
  generality, $L/K$ has degree $n$, and $\xi$ is the primitive
  $n$th root of unity $\zeta$ chosen in the construction of the
  cohomology class $C^{a,b}$.

  We use $\gamma$ to construct a matrix $A$ as in
  Lemma~\ref{thm:lemma-let-b}. Let $G$ be the diagonal matrix with
  $\gamma\gamma^\sigma\dots\gamma^{\sigma^{i-1}}$ in the $(i,i)$
  position. Let $x_1, \dots , x_n$ be any $K$-basis for $L$, and let $X$
  be the matrix with $x_j^{\sigma^{i-1}}$ in the $(i,j)$ position. Let
  $A= GX$, so that
  \begin{equation*}
    A_{i,j} =
    \gamma\gamma^{\sigma}\dots\gamma^{\sigma^{i-1}}x_j^{\sigma^{i-1}},
  \end{equation*}
  and thus
  \begin{equation*}
    \gamma A_{i,j}^\sigma =
    \begin{cases}
      A_{i+1, j} & i < n\\
      b A_{1, j} & i = n.
    \end{cases}
  \end{equation*}
  From this, we see that $\gamma A^\sigma = M_b A$, that is,
  $$A^\sigma A^{-1} \equiv C_\sigma^{a,b} \pmod{\McOmega^\times}.$$
  This induces an equality of cocycles, as $\sigma$ generates
  $\Gal(L/K)$.
  Since
  \begin{equation*}
    \gamma\gamma^\sigma \dots \gamma^{\sigma^{i-1}}\McO{L} =
    \mathfrak{c} \mathfrak{c}^{-\sigma^i} \mathfrak{b}^i,
  \end{equation*}
  we have
  $$
    \det(G)\McO{L}
    = \mathfrak{c}^n(\Norm_{L/K} \mathfrak{c})^{-1}
    \mathfrak{b}^{n(n+1)/2}.
  $$
  Thus, reducing exponents modulo $n$ and using
  Lemmas~\ref{cupandazumaya} and~\ref{thm:lemma-let-b}, we find
  \begin{equation}\label{eq:mcpairing}
    \Mcpairing{a}{-ab} =\delta C^{a,b} \tensor \zeta =
    ([\Norm_{L/K}\mathfrak{c}]+ \frac{n}{2}[\mathfrak{b}]
    - [\det(X)\McO{K}]) \tensor \zeta,
  \end{equation}
  where, by a slight abuse of notation, $\det(X)\mathcal{O}_{K,S}$
  denotes the ideal of $\mathcal{O}_{K,S}$ whose extension to
  $\mathcal{O}_S$ is $\det(X)\mathcal{O}_S$, which exists because
  $\det(X)^\sigma = \pm \det(X)$.
  As a special case, we take $b=1$, so $\mathfrak{c}=\mathfrak{b}=1$
  and
  $$
    \Mcpairing{a}{-a} = -[\det(X)\McO{K}] \tensor \zeta.
  $$
  Subtracting this from~\eqref{eq:mcpairing},
  we obtain the theorem.
\end{proof}

We have the following corollaries (which are also easy to prove
directly).

\begin{corollary}\label{PairingisZero}
  If $a, b \in \McU{K}$, and if $b$ is the norm of an element of
  $\McU{L}$, then $\Mcpairing{a}{b} = 1$.
\end{corollary}

\begin{corollary}\label{MilnorK}
  If $a\in \McU{K}$ is such that $1-a \in \McU{K}$, then
  $\Mcpairing{a}{1-a} = 1$.
\end{corollary}

\section{The cup product and $K$-theory} \label{K-theory}

Corollary~\ref{MilnorK} may be rephrased in terms of $K$-theory.
Recall the definition of the Milnor $K_2$-group of a commutative
ring $R$:
$$
  K_2^M(R) = (R^{\times} \otimes R^{\times})/ \langle a \otimes (1-a)
  \colon a,1-a \in R^{\times} \rangle.
$$
Then Corollary~\ref{MilnorK} says that the restriction of the cup
product to the $S$-units induces a map
$$
  u \colon K_2^M(\McO{K})/n \to H^2(G_{K,S},\mu_{n}^{\otimes 2}).
$$
On the other hand, since $\mu_{n} \subset K$, the exact sequence
of Tate \cite[Theorem 6.2]{tate:1976} and \eqref{H2isAmodp} yield
a (noncanonical) isomorphism
$$
  c \colon K_2(\McO{K})/n \xrightarrow{\sim}
  H^2(G_{K,S},\mu_{n}^{\otimes 2})
$$
(see \cite{soule:1979}, \cite{dwyer-friedlander:1985} and
\cite{keune:1989} for generalizations).  In \cite{soule:1979}, a
particular choice of the map $c$ is described as a Chern class map
(for $n$ a power of a prime $p$).  The two versions of $K_{2}$ are
related by a map
\begin{equation} \label{Kcup}
  \kappa \colon K_2^M(\McO{K}) \to K_2(\McO{K}),
\end{equation}
which is constructed as follows.  First, we may identify $K_2^M(K)$
with $K_2(K)$ by a classical result of Matsumoto.  The group
$K_2(\McO{K})$ may then be defined via the exact localization sequence
\begin{equation} \label{localseq}
  0 \to K_2(\McO{K}) \to K_2(K) \xrightarrow{t}
  \bigoplus_{\q \notin S} k_{\q}^{\times} \to 0,
\end{equation}
where $k_{\q}$ denotes the residue field of $K$ at ${\q}$, and the
map $t$ is given by tame symbols.  Since two $S$-units pair
trivially under the tame symbols, the sequence \eqref{localseq}
yields the map $\kappa$ of \eqref{Kcup}.  Let $\kappa_{n}$ denote
the map induced by $\kappa$ on $K$-groups modulo $n$.

\begin{proposition} \label{Krelation}
  The maps $u$, $c$, and $\kappa$ are related by the commutative
  diagram
\[
\xymatrix{ K_2^M(\McO{K})/n \ar[r]^{\kappa_{n}} \ar[rd]^{u} &
  K_2(\McO{K})/n
  \ar[d]^{-c}_{\iso} \\
  & H^2(G_{K,S},\mu_{n}^{\otimes 2}).  }
\]
\end{proposition}

\begin{proof}
  This is a special case of \cite[Theorem 1]{soule:1979} for $n$ a
  prime power, and the general case follows easily.
\end{proof}

Thus, the question of whether the pairing $\blankpairing$ is
surjective on $S$-units is the question of surjectivity of
$\kappa_{n}$.  In Section~\ref{cyclpunits}, we make a conjecture
on the surjectivity of $\kappa_{p}$
(Conjecture~\ref{k2conjecture}) for the ring of $p$-integers in
$\QQ(\mu_p)$. In order for this conjecture to hold, $\kappa_{p}$
will often have to be injective as well.  We now describe a
necessary, general condition for this injectivity.

The Steinberg symbol $\{a, b\}$, for $a, b \in R^{\times}$, is defined
to be the image of $a \otimes b$ in $K_2^M(R)$.  For any field $F$,
the group $K_2^M(F)$ has the property that the Steinberg symbols are
antisymmetric.  This need not be true in an arbitrary ring $R$.  On
the other hand, in order that $\kappa_{n}$ be injective, it is
necessary that antisymmetry hold in $K_2^M(\McO{K})/n$, since
$K_2(\McO{K}) \subset K_2^M(K)$.  We now present one sufficient
condition for this antisymmetry.

\begin{lemma} \label{antisymmetry}
  Fix $\varepsilon \in \ZZ$ relatively prime to $n$.  Assume that
  $R^{\times}/R^{\times n}$ has a generating set with a set of
  representatives $S \subseteq R^{\times}$ such that
  \begin{equation} \label{unitproperty}
     \{ 1-st^{\varepsilon} \colon s \in S,\ t \in S \cup \{1\},\
     s \neq t \}  \subset R^{\times}.
  \end{equation}
  Then $\{a,b\} + \{b,a\} \equiv 0 \pmod{n}$ for any $a, b \in
  R^{\times}$.
\end{lemma}

\begin{proof}
  Let $s \in S$ and $t \in S^{\varepsilon} \cup \{1\}$ with $s \neq
  t$, and set $x = st$.  By definition of $K_2^M(R)$, we have
  \[
  \left\{x,-x\right\} = -\left\{x,-\frac{1}{x}\right\} =
  -\left\{x,\frac{x-1}{x}\right\} =
  \left\{\frac{1}{x},\frac{x-1}{x}\right\} = 0.
  \]
  Then
  \[
  \{s,t\} + \{t,s\} = \{s,-st\} + \{t,-st\} = \{st,-st\} = 0.
  \]
  In general, if $a = x^{n}\prod_{i=1}^{M} a_i$ and $b =
  y^{n}\prod_{j=1}^{N}b_j$ with $a_i, b_j \in S \cup S^{-1}$ and
  $x,y \in R^{\times},$ then
  \[
  \{a,b\} + \{b,a\} \equiv \sum_{i=1}^M\sum_{j=1}^N (\{a_i,b_j\} +
  \{b_j,a_i\}) \equiv 0 \pmod{n}.
  \]
\end{proof}

\section{Pairing with a $p$th root of unity}
\label{IwasawaLfunction}

We now focus on the case $n = p$, an odd prime, $K = \QQ(\mu_p)$
and $S = \{(1-\zeta)\}$.  Here, and for the remainder of the
paper, we fix a choice $\zeta$ of a primitive $p$th root of unity.
In this section, we will show that the cup product we are
considering is nontrivial for many $p$.  We assume Vandiver's
conjecture holds at $p$, i.e., that $p$ does not divide the class
number of $\QQ(\zeta+\zeta^{-1})$. Let $\mathcal{C}$ denote the
group of cyclotomic $p$-units.

\begin{lemma}\label{Cyclotomic}
  The symbol $\{\zeta,x\} \in K_2^M(\McO{K})$ is zero for all $x \in
  \McU{K}$.  In particular, $\Mczeta$ pairs trivially with $\McU{K}$
under the cup product.
\end{lemma}

\begin{proof}
  This follows immediately from the properties of the symbol, since
  the elements $1-\zeta^i$ with $1 \le i \le p-1$ generate
  $\mathcal{C}$, which, by a well-known consequence of
  Vandiver's conjecture,
 has index prime to $p$ in $\McU{K}$.
\end{proof}

Let $A_K$ denote the $p$-part of the class group $\clnos{K}$ of $K$.
(Note that $C_K = C_{K,S}$.)  We have an exact sequence
$$
0 \to \McU{K}/\McO{K}^{\times p} \to H^{1}(G_{K,S}, \mu_{p}) \to
A_K[p] \to 0,
$$
in which the map on the right is induced by $a \mapsto {\mathfrak
  a}$, where $a\McO{K} = \mathfrak{a}^p$. By Lemma~\ref{Cyclotomic},
the map $a \mapsto \Mcpairing{\zeta}{a}$ factors through a map
\[
\psi \colon A_K[p] \to A_K \tensor \mu_{p}.\label{MaponA}
\]
Let $K_{\infty}/K$ be the cyclotomic $\ZZ_{p}$-extension, and let
$A_{\infty}$ be the inverse limit of the $p$-parts of the ideal
class groups under norm maps up the cyclotomic tower.  The abelian
group $A_{\infty}$ breaks up into eigenspaces
$A_{\infty}(\omega^{i})$, on which $\Gal(K/\QQ)$ acts by the $i$th
power of the Teichm\"uller character $\omega \colon \Delta \to
\zp^{\times}$.  Let $\Lambda$ be the Iwasawa algebra.  Fix a
topological generator $\gamma$ of $\Gal(K_{\infty}/K)$ such that
$\gamma$ acts on $\mu_{p^2}$ by raising to the ($1+p$)th power,
and let $T = \gamma-1$ be the corresponding variable in the
Iwasawa algebra, so $\Lambda \iso \ZZ_{p}[[T]]$. Since $p$
satisfies Vandiver's conjecture, $A_{\infty}(\omega^{i}) \iso
\Lambda/(f_{i})$, where $f_{i}$ is a characteristic power series
\cite{washington:1997}.

We consider an $i$ for which $A_K(\omega^i)$ is nontrivial.  We
often abuse notation by using the same symbol to denote both a
(fractional) ideal and its ideal class.

\begin{proposition} \label{zetapairing}
  Let $\mathfrak{a} \in A_K[p](\omega^i)$, and choose $\mathfrak{a}_0
  \in A_K(\omega^i)$ with $\mathfrak{a}_0^{f_i(0)/p} = \mathfrak{a}$.
  Then
  $$
  \psi(\mathfrak{a}) = \mathfrak{a}_0^{-f_{i}'(0)} \tensor \Mczeta.
  $$
\end{proposition}

\begin{proof}
  Let $L = K(\mu_{p^2})$, set $f = f_i$ and $\Norm = \Norm_{L/K}$, and
  let ${\mathfrak{a}}_{1}$ be an ideal of $\McO{L}$ with norm
  ${\mathfrak{a}}_{0}$.  The ideal ${\mathfrak{a}}_{1}^{f(T)}$ is
  principal, generated by an element $y$, and
  \begin{equation}
   (\Norm y)\McO{K} = \Norm{\mathfrak{a}}_{1}^{f(0)} =
   {\mathfrak{a}}^{p} = a\McO{K}
  \label{eq:1}
  \end{equation}
  for some $a \in D_K$.  Since $\zeta$ pairs trivially with $p$-units
  by Lemma~\ref{Cyclotomic}, we may assume that $y$ is the element we
  use to calculate $\Mcpairing{\Mczeta}{a}$ in
  Theorem~\ref{PairingFormula}.  In particular, we have that
  \begin{equation}
        y\McO{L} = {\mathfrak{a}}_{1}^{f(T)} =
        {\mathfrak{a}}{\mathfrak{b}}^{T}\label{eq:10}
  \end{equation}
  for ${\mathfrak{b}}$ such that
  $$
  \Mcpairing{\Mczeta}{a} = \Norm{\mathfrak{b}}^{-1} \otimes
  \Mczeta.
  $$

  Let us make the identification
  $$
  A_L(\omega^i) \iso \ZZ_{p}[[T]]/(f(T),(1+T)^{p}-1)
  $$
  by mapping ${\mathfrak{a}}_{1}$ to $1$.  Then, using~\eqref{eq:1}
  and \eqref{eq:10} and that $\Norm = ((1+T)^p-1)/T$, we see that
  ${\mathfrak{b}}$ may be identified with
  $$
  \frac{1}{T}\left(f(T) - \frac{(1+T)^p - 1}{pT}f(0)\right) \equiv
  f'(0) \pmod{(p,T)},
  $$
  and hence
  $$
  \Norm{\mathfrak{b}} \equiv {{\mathfrak{a}}_0}^{f'(0)}
  \pmod{pA_K}.
  $$
\end{proof}

Thus if $f_i'(0)$ is nonzero modulo $p$, which is equivalent to
saying that the $\lambda$-invariant of $A_{\infty}(\omega^i)$ is
$1$, then the pairing $\blankpairing$ is nontrivial.  This occurs
for all irregular primes less than $12,\!000,\!000$
\cite{buhler-crandall-ernvall-metsankyla-shokrollahi-2001}.
(We remark that for many purposes in this article, such as
Proposition \ref{zetapairing}, the cyclicity of $A_K(\omega^i)$
would be a sufficient assumption \cite{kurihara:1993}, but we
are content to assume the stronger condition of Vandiver's
conjecture.)

\section{Restriction of the pairing to the cyclotomic $p$-units}
\label{cyclpunits} We continue to assume $K = \QQ(\mu_p)$, $S = \{
(1-\zeta) \}$ and $n = p$, an odd prime satisfying Vandiver's
conjecture.  We consider the restriction of the pairing
$\blankpairing$ to the subgroup ${\mathcal
  C}$ of cyclotomic $p$-units.

A pair $(p,r)$, with $r$ even and $2 \le r \le p-3$, is called
irregular if $p$ divides the numerator of the $r$th Bernoulli number
$B_r$, that is, if the eigenspace $A_K(\omega^{p-r})$ is nontrivial.
Given such a pair, we choose an isomorphism
\begin{equation} \label{isochoice}
        A_K(\omega^{p-r}) \tensor \mu_{p} \iso \ZZ/p\ZZ(2-r).
\end{equation}
(We shall view the underlying group structure of $A(i)$ for a
given $G_{K,S}$-module $A$ and twist $i \in \ZZ$ as being
canonically identified with that of $A$.) Let $\Delta =
\Gal(K/\QQ)$, and consider the $\Delta$-equivariant pairing
$$
\blankeigpair{r}: D_K \times D_K \to \ZZ/p\ZZ(2-r)
$$
arising via \eqref{isochoice} from composition of $\blankpairing$
with projection onto $A_K(\omega^{p-r}) \tensor \mu_{p}$.

We consider the restriction of $\blankeigpair{r}$ to $\mathcal{C}
\times \mathcal{C}$.  In fact, since $\mathcal{C} = \langle -\zeta
\rangle \oplus \mathcal{C}^{+}$, we can by Lemma~\ref{Cyclotomic}
focus attention on the restriction of the pairing to elements of
$\mathcal{C}^{+}$.  Eigenspace considerations put some restrictions on
the elements that can pair nontrivially.  For any integer $i$,
consider the usual idempotent
$$
\epsilon_{i} = \frac{1}{p-1}\sum_{\sigma \in \Delta}
\omega(\sigma)^{-i}\sigma,
$$
and choose $\eta_i \in \mathcal{C}$ with
$$
\eta_i \equiv (1-\zeta)^{\epsilon_{p-i}} \pmod{\mathcal{C}^p}.
$$
Then $\mathcal{C}^+$ is generated by $\{\eta_{i}:\mbox{$i$ odd}, 1
\le i \le p-2\}$. Since $\Delta$ acts on
$\Mcpairing{\eta_{i}}{\eta_{j}}$ via $\omega^{2-i-j}$, we have
\begin{equation*}
   \mbox{$\eigpair{\eta_{i}}{\eta_{j}}{r} = 0$ if
   $i+j \not\equiv r\pmod{p-1}$}.
\end{equation*}

Let
\begin{equation*}
        e_{i,r} = \eigpair{\eta_{i}}{\eta_{r-i}}{r},
        \quad \text{$i$ odd, $1 \le i \le p-2$}.
\end{equation*}
We now impose many linear relations on the $e_{i,r}$ and, in the
process, get bounds on the order of $K_2^M(\McO{K})/p$.
Lemma~\ref{Cyclotomic} and Corollary~\ref{MilnorK} imply that if
$\eta$ is a cyclotomic unit and $1-\eta = \Mczeta^{j}\eta'$ for some
$j$ and cyclotomic unit $\eta'$, then $\{\eta,\eta'\} = 0$ in
$K_2^M(\McO{K})$, and hence $(\eta,\eta')_S = 1$.  Applying this
observation to
$$
\eta = \rho_{a} = \sum_{j=0}^{a-1} (-\zeta)^j,
$$
we get
$$
\{\rho_{a},\rho_{a-1}\} = 0, \qquad 3 \le a \le p-1.
$$
Note that, if $a$ is even,
$$
\rho_{a} = \frac{1-\zeta^{a}}{1+\zeta} =
\frac{(1-\zeta^{a})(1-\zeta)}{1-\zeta^{2}}, \quad \rho_{a-1} =
\frac{1+\zeta^{a-1}}{1+\zeta} =
\frac{(1-\zeta^{2a-2})(1-\zeta)}{(1-\zeta^{a-1})(1-\zeta^{2})}.
$$
If $\sigma \in \Delta$ satisfies $\sigma\zeta = \zeta^a$, then
$\omega(\sigma) \equiv a \pmod{p}$.
Thus
$$
(1-\Mczeta^{a})^{\epsilon_{p-i}} \equiv \eta_{i}^{a^{p-i}}
\pmod{\mathcal{C}^p}.
$$
Hence, the $e_{i,r}$ must be solutions $x_i = e_{i,r}$ over
$\ZZ/p\ZZ$ to
\begin{equation} \label{Relations}
     \sum_{\substack{\text{$i$ odd} \\ 1 \le i \le p-2}}
     (1+a^{p-i}-2^{p-i})(1-2^{p-r+i})(1-(a-1)^{p-r+i})x_i = 0
\end{equation}
for every even $a$ with $4 \le a \le p-1$.  (These relations also hold
for odd $a$ with $3 \le a \le p-2$, but we will not use those.)

\begin{theorem} \label{uniquepairing}
  For all irregular pairs $(p, r)$ with $p < \maxp$, there exists a
  nontrivial, Galois equivariant, skew-symmetric pairing
  $$\blankeigpair{} \colon \mathcal{C} \times \mathcal{C} \to
  \ZZ/p\ZZ(2-r)$$
  satisfying \eqref{Relations} with $x_i =
  \eigpair{\eta_i}{\eta_{r-i}}{}$.  Furthermore, these properties
  uniquely define the pairing up to a scalar multiple.
\end{theorem}

\begin{proof}[Sketch of Proof]
  The relations \eqref{Relations} and the antisymmetry relations, $x_i
  + x_{r-i} = 0$, put constraints on possible values of the pairing.
  We used a computer to calculate the nullspace of the matrix of
  coefficients in these relations.
\end{proof}

We have computed the pairing of Theorem~\ref{uniquepairing} for
all $(p,r)$ with $p < \maxp$.  A table of the pairings and Magma
routines that perform the computation are available at {\tt
www.math.harvard.edu/{\urltilde}sharifi} and {\tt
www.math.arizona.edu/\urltilde{}wmc}.  The pairing
$\blankeigpair{r}$ must be a (possibly zero) scalar multiple of
the computed pairing.

\begin{corollary} \label{k2corollary}
  For all irregular pairs $(p,r)$ with $p < \maxp$, one has
  \[ |(K_2^M(\McO{K})/p)(\omega^{2-r})| \le p. \]
\end{corollary}

\begin{proof}
  The only point at issue here is whether the symbols in
  $K_2^M(\McO{K})/p$ satisfy the skew-symmetry that was used in the
  proof of Theorem~\ref{uniquepairing}, as the other relations used in
  Theorem~\ref{uniquepairing} arise directly from relations in
  $K_2^M(\McO{K})$ and eigenspace considerations.  Since $\mathcal{C}$
  has index prime to $p$ in $\McU{K}$, we need only remark that the
  set $S$ of generators of $\mathcal{C}$ given by $1-\zeta^i$ with $1
  \le i \le p-1$ has the property \eqref{unitproperty} for
  $\varepsilon = -1$.  Then, by Lemma~\ref{antisymmetry}, the image of
  $a \otimes b + b \otimes a$ is trivial in $K_2^M(\McO{K}) \otimes
  \zp$ for any $a,b \in \McU{K}$.
\end{proof}

In fact, by performing modulo $p^2$ the same computations as in the
proof of Theorem~\ref{uniquepairing}, we have verified that
\[
|(K_2^M(\McO{K}) \otimes \zp)(\omega^{2-r})| \le p
\]
for all irregular pairs with $p < 3000$.  The authors suspect that the
uniqueness in Theorem~\ref{uniquepairing} fails for the irregular pair
given by $p = 89209$ and $r = (p+3)/2$, as the values $x_1 =
-x_{(p+1)/2} = 1$ and $x_i = 0$ for $i \neq 1,(p+1)/2$ provide a
solution to the equations in this case, since $2^{p-r+1}\equiv
1\pmod{p}$.  Corollary~\ref{k2corollary} may or may not still hold for
this $p$, since we did not use all of the defining relations of
$K_2^M(\McO{K})$ in its proof.  We do, however, conjecture that the
pairing $\blankeigpair{r}$ is nontrivial.  Since the eigenspaces of
$K_2(\McO{K})/p$ corresponding to regular pairs are trivial, this may
be rephrased as follows.

\begin{conjecture} \label{k2conjecture}
  Let $p$ be an odd prime satisfying Vandiver's conjecture.  Let $K =
  \QQ(\mu_p)$ and let $\McO{K}$ denote the ring of $p$-integers in
  $K$.  The natural map
  $$
  K_2^M(\McO{K}) \otimes \zp \to K_2(\McO{K}) \otimes \zp
  $$
  is surjective.
\end{conjecture}

In Section~\ref{AnotherFormula}, we verify this conjecture in the
case $p=37$.  In general, we have constructed many other relations
similar to those of \eqref{Relations}, and their solutions are
consistent with the values of the calculated pairing of
Theorem~\ref{uniquepairing} in those cases that we have tested. As
the example with $p = 89209$ illustrates (or consider the regular
pair $(73,38)$), we have no convincing evidence regarding whether
or not the map in Conjecture~\ref{k2conjecture} is always
injective.

As further circumstantial evidence for the conjecture, we note
that the nontrivial pairing of Theorem~\ref{uniquepairing} has a
property not obviously encoded in the relations above, namely that
$x_{p-r} = 0$.
The cup product pairing $\blankeigpair{r}$ must itself satisfy this
relation, since $\eta_{p-r}$ provides a Kummer generator for the
unramified extension of $K$ whose Galois group corresponds to
$(A_K/pA_K)(\omega^{p-r})$ with respect to the Artin map.  Thus the
norm of an ideal from this extension always has trivial projection to
the $\omega^{p-r}$ eigenspace of $A_K/pA_K$, and hence the pairing
must be trivial by Theorem~\ref{PairingFormula}.

\section{Relationship with the ideal class group} \label{idealclass}

Let $K$ be a number field and $S$ a set of primes of $K$ that
contains all real places of $K$. For a finite extension $F$ of
$K$, we denote by $\McI{F}$, $\McP{F}$, and $\cl{F}$ the ideals,
principal ideals, and ideal class group, respectively, of the
$S$-integers $\McO{F}$, and by $\McH{F}$ the maximal unramified
abelian extension of $F$ in which all primes above $S$ split
completely.

First we review some genus theory. Let $L/K$ be a cyclic
extension, unramified outside $S$, with Galois group $G$ generated
by
an element
$\sigma$.
\begin{lemma}\label{genustheory}
  The norm map $\Norm_{L/K}: \cl{L} \to \cl{K}$ induces a map
  \begin{equation*}
    \cl{L}/(\sigma-1)\cl{L} \to \cl{K},
  \end{equation*}
  which is a surjection if and only if $L \cap \McH{K} = K$, and is an
  injection if there is at most one prime in $S$ that does not split
  completely in $L/K$.
\end{lemma}

\begin{proof}
  By class field theory, $\cl{L} = \Gal(\McH{L}/L)$ and $\cl{K} =
  \Gal(\McH{K}/K)$.  With these identifications, the norm map is
  restriction to $\McH{K}$.  Thus the image of the norm map is
  $\Gal(\McH{K}/(L \cap \McH{K}))$, which immediately implies the first
  assertion.
  The kernel of the norm map is generated by the
  commutator subgroup and
  the intersection between
  $\Gal(\McH{L}/L)$ and the subgroup of $\Gal(\McH{L}/K)$ generated by
  decomposition groups of primes in $S$.
  If there is exactly one nontrivial decomposition group, then its
  contribution is trivial, since it maps injectively to $\Gal(L/K)$.
  Furthermore, since $G$ is
  cyclic the commutator subgroup of $\Gal(\McH{L}/K)$ is
  $[\sigma,\Gal(\McH{L}/L)] = (\sigma-1)\cl{L}$.
\end{proof}

Compare the following proposition with
Theorem~\ref{PairingFormula}.

\begin{proposition}\label{pairingandcapitulation}
  There is an isomorphism
  $$
  \cl{L}^G/\phi(\cl{K}) \xrightarrow{\sim} (\McU{K} \cap
  \Norm_{L/K} L^{\times})/ \Norm_{L/K} \McU{L}
  $$
  given by taking the class of an ideal $\mathfrak{a}$ to an
  element $b \in \McU{K}$ such that $b = \Norm_{L/K} y$ with $y\McO{L}
  = \mathfrak{a}^{1-\sigma}$.  If $\McH{K} \cap L = K$ and there is at
  most one prime in $S$ that does not split completely in $L/K$,
  the order of these groups is equal to the number of ideal classes in
  $\cl{K}$ that define trivial classes in $\cl{L}$.
\end{proposition}

\begin{proof}
  Consider the exact sequence
  $$
  0 \to \McP{L} \to \McI{L} \to \cl{L} \to 0.
  $$
  Since $L/K$ is unramified outside $S$, $\McI{L}$ is a direct sum
  of induced modules, so $\hat{H}^{1}(G,\McI{L}) = 0$. Thus we have a
  surjection
\begin{equation*}
  \hat{H}^{0}(G, \cl{L}) \to \hat{H}^{-1}(G, \McP{L}),
\end{equation*}
which takes the class of an ideal $\mathfrak{a}$ to
$({1-\sigma})\mathfrak{a}$.  Since the map $\McP{L}^G \to
\McI{L}^G = \McI{K}$ has cokernel equal to the image of $\cl{K}$
under the obvious map $\phi \colon \cl{K} \to \cl{L}$, we
therefore have an isomorphism
$$
\cl{L}^G/\phi(\cl{K}) \iso \hat{H}^{-1}(G,\McP{L}).
$$
Furthermore, using the exact sequence
$$
0 \to \McU{L} \to L^{\times} \to \McP{L} \to 0
$$
and the triviality of $\hat H^{-1}(G, L^{\times})$, we obtain
$$
\hat H^{-1}(G, \McP{L}) \iso (\McU{K}\cap
\Norm_{L/K}L^{\times})/\Norm_{L/K}\McU{L},
$$
induced by $\Norm_{L/K}$ on a generator of a representative
principal ideal.  This proves the first statement in the
proposition.

Now, suppose that $\McH{K} \cap L = K$ and that there is at most
one prime in $S$ that does not split completely in $L/K$.  Then it
follows from Lemma~\ref{genustheory} that $\phi(\cl{K}) =
N_G(\cl{L})$, where $N_G$ is the norm element in the group ring
$\ZZ[G]$. Thus
$$\cl{L}^G/\phi(\cl{K}) = \hat H^0(G, \cl{L}).$$
Furthermore, since $\cl{L}$ is finite, we have
$$
|\hat H^{0}(G, \cl{L})| = |\hat H^{-1}(G, \cl{L})|.
$$
Finally, using Lemma~\ref{genustheory} again, we see that
$$
\hat H^{-1}(G, \cl{L}) \iso \mathrm{ker}(\cl{K} \to \cl{L}).
$$
This proves the second assertion of the proposition.
\end{proof}

Now we return to the situation of Section~\ref{cupproduct}, fixing
$n$ and letting $K$ contain the $n$th roots of unity and $S$ all
primes above $n$ and real archimedean places.  Let $a \in D_K$ and
$\alpha \in K_S^{\times}$ with $\alpha^{n} = a$ and take $L =
K(\alpha)$.  We let $d = [L:K]$.

\begin{proposition} \label{nonzeropairing}
  Suppose $a \in D_K$ is such that $\McH{K} \cap L = K$ and there is
  at most one prime of $S$ that does not split completely in $L/K$.  Then
  the map $(n/d)\Norm_{L/K}: \cl{L} \rightarrow \cl{K}$ induces an
  isomorphism
    \begin{equation} \label{isomorphicgroup}
        \cl{L}^G/((d,\sigma-1)\cl{L})^G \otimes \mu_n \xrightarrow{\sim}
        \Mcpairing{a}{\McU{K}}\cap (\cl{K}\tensor \mu_n),
    \end{equation}
    the intersection being taken in $H^2(G_{K,S}, \mu_n)\tensor
    \mu_n$.
\end{proposition}

\begin{proof}
  Given $\mathfrak{a} \in\cl{L}^G$, we can find $b \in \McU{K} \cap
  N_{L/K}L^\times$ associated with $\mathfrak{a}$ by the isomorphism in
  Proposition~\ref{pairingandcapitulation}.  Since $b$ is a global
  norm it is a local norm everywhere, and Theorem~\ref{PairingFormula}
  applies.  Conversely, given $b \in \McU{K}$ such that $(a,b)_v = 1$
  for all $v\in S$, Theorem~\ref{PairingFormula} supplies an ideal
  $\mathfrak{a}$ with class in $\cl{L}^G$ such that $\Mcpairing{a}{b} =
  \Norm_{L/K} \mathfrak{a} \otimes \xi$, with $\xi$ a fixed generator
  of $\mu_d$.  Therefore, the map in \eqref{isomorphicgroup} is
  surjective.  By Lemma~\ref{genustheory}, the kernel of $\Norm_{L/K}$
  on $\cl{L}$ is $(\sigma-1)\cl{L}$ and $\Norm_{L/K}$ is surjective.
  Thus, the kernel of the map $\cl{L} \to \cl{K} \otimes \mu_n$ given
  by $\mathfrak{a} \mapsto N_{L/K} \mathfrak{a} \otimes \xi$ is
  $(d,\sigma-1)\cl{L}$.
\end{proof}

This has, for instance, the following corollary in the case $n =
p$, a prime number. Let $\pcl{F}$ denote the $p$-part of the
$S$-class group of $F$ for $F/K$ finite.

\begin{corollary} \label{simplecase}
  Assume that $|\pcl{K}| = p$ and $S$ consists of a single
  $($unique$)$
  prime above $p$.  Let $a \in D_K$ be such that $[L:K] = p$. Then
  $\Mcpairing{a}{\McU{K}} \neq 0$ if and only if $|\pcl{L}| = p$.
\end{corollary}

\begin{proof}
  If $L \subseteq \McH{K}$, then since $\pcl{L}$ is the commutator
  subgroup of a $p$-group that has maximal abelian quotient $\pcl{K}
  \iso \ZZ/p\ZZ$, we must have $\pcl{L} = 0$.  Hence, we may assume
  $\McH{K} \cap L = K$.

  By Lemma~\ref{genustheory}, the order of the quotient
  $\pcl{L}/(\sigma-1)\pcl{L}$ is $p$.  By the assumption on $S$, the
  image of the pairing is contained in $\cl{K}\tensor \mu_p$.  Thus,
  using Proposition~\ref{nonzeropairing}, we see that $\pcl{L}^G$
  surjects onto the above quotient if and only if
  $\Mcpairing{a}{\McU{K}}$ is nonzero.  On the other hand, $\pcl{L}^G$
  surjects onto the quotient if and only if $(\sigma-1)\pcl{L} = 0$.
\end{proof}

Again, let us consider the case $n = p$ odd, $K = \QQ(\mu_p)$ and
$S = \{ (1-\zeta) \}$.  Assume Vandiver's conjecture, and let
$(p,r)$ be an irregular pair.  Recall the pairing
$\blankeigpair{r}$ of Section~\ref{cyclpunits}.

\begin{lemma} \label{pairingimage}
  The image of the pairing $\blankeigpair{r}$ is
  $\eigpair{1-\zeta}{\mathcal{C}^{+}}{r}$.
\end{lemma}

\begin{proof}
  The image of the pairing is generated by
  $\eigpair{\eta_i}{\eta_{r-i}}{r}$ for all odd $i$. Since
  $\eigpair{\eta_i}{\eta_j}{r} = 0$ for $j \not\equiv r-i \mod p-1$,
  we have
  $$
  \eigpair{\eta_i}{\eta_{r-i}}{r} =
  \eigpair{\eta_i}{\prod_{j=0}^{p-2}(1-\zeta)^{\epsilon_j}}{r} =
  \eigpair{\eta_i}{1-\zeta}{r}.
  $$
\end{proof}

Proposition~\ref{nonzeropairing} allows us to conclude the
following.

\begin{corollary} \label{specialcase}
  Let $\alpha^p = 1-\zeta$ and $L = K(\alpha)$.  Then the image of the
  pairing $\blankeigpair{r}$ is isomorphic to the
  $\omega^{p-r}$-eigenspace of $A_L^G/((p,\sigma-1)A_L)^G$, where
  $A_L$ denotes the $p$-part of the class group of $L$.
\end{corollary}

\section{Relationship with local pairings}\label{AnotherFormula}

We restrict ourselves to the case $n = p$ and $K$ containing
$\mu_p$. We assume that $S$ consists of a single, unique prime of
$K$ above $p$ and, if $p=2$, that $K$ has no real places. In
Theorem~\ref{HCFpairformula}, we will derive a formula for (a
projection of) the pairing $(a,b)_S$ as a norm residue symbol in a
certain unramified cyclic extension $L$ of $K$ of degree $p$.  The
formula is similar to that of Theorem~\ref{PairingFormula} in that
its applicability amounts to the determination of an element $c
\in L^{\times}$ such that $c^{\sigma-1}b \in L^{\times p}$, where
$\sigma$ generates $\Gal(L/K)$ (for a different field $L$). In
this case, however, one must determine an embedding of $c$ in the
multiplicative group modulo $p$th powers of the completion at $L$
at a prime above $p$, as opposed to determining the class modulo
$p$ of an ideal of $\McO{K}$ that $c$ generates in $L$, again up
to a $p$th power.

Let $L/K$ be an unramified cyclic extension of degree $p$, and set
$G = \Gal(L/K)$.  Consider the following commutative diagram, in
which \eqref{H2isAmodp} has been used to obtain the top row and in
which we have identified $H^2(G_{K,S}, \mu_p)$ with
$\cl{K}/p\cl{K}$ in the bottom row.
\begin{equation*}
        \xymatrix{0 \ar[r] & \cl{L}/p\cl{L} \ar[r]\ar[d]^{\Norm_{L/K}} & H^{2}(G_{L,S}, \mu_{p})
        \ar[r]\ar[d]^{\cor} \ar[r] &
        \bigoplus^0_{\p|p} H^{2}(L_{\p},\mu_{p})\ar[d]^{f}\ar[r] &0\\
        0 \ar[r] & \Norm_{L/K}\cl{L}/p\cl{K} \ar[r] &
        \cl{K}/p\cl{K}\ar[r] & \cl{K}/\Norm_{L/K}\cl{L} \ar[r]&0
        }
\end{equation*}
Here, the superscript 0 on the direct sum indicates the kernel of
the map
\begin{equation*}
  \sum_{\p|p} \inv_{\p}: \bigoplus_{\p|p}
  H^{2}(L_{\p},\mu_{p}) \rightarrow
  \mathbb{Z}/p\mathbb{Z}
\end{equation*}
(under the obvious identification $\frac{1}{p}\ZZ/\ZZ \iso
\ZZ/p\ZZ$). By class field theory, there is a natural isomorphism
$\cl{K}/\Norm_{L/K}\cl{L} \iso G$. Let $\sigma$ be a generator of
$G$. Choose a prime of $L$ above $p$, say $\p_0$, and for any
other, let $k_{\p} \in \mathbb{Z}/p\mathbb{Z}$ be such that $\p =
\sigma^{-k_{\p}} \p_0$. We need the following explicit description
of the map $f$.

\begin{lemma}\label{formulaforf}
  The map
        \begin{equation*}
                f: \bigoplus^0_{\p|p} H^{2}(L_{\p},\mu_{p}) \to
        \cl{K}/\Norm_{L/K}\cl{L}
        \end{equation*}
        is given by
        $$
        f(c) = \sum_{\p\mid p}k_{\p} \inv_{\p}(c) \cdot
        \mathfrak{c}
        $$
        for some ideal class $\mathfrak{c}$ generating
        $\cl{K}/N_{L/K}\cl{L}$.
\end{lemma}

\begin{proof}
  Corestriction is equivariant with respect to $G$, and hence so is
  the map $f$.  Furthermore, $\bigoplus^0 H^{2}(L_{\p}, \mu_{p})$ is a
  cyclic $\ZZ_{p}[G]$-module, and
  $$
    \cl{K}/\Norm_{L/K}\cl{L} \iso \mathbb{Z}/p\mathbb{Z},
  $$
  so there is
  only one non-zero map up to scalar multiple.
  Since
  $G_{K,S}$ has cohomological dimension at most 2 \cite[Proposition
  8.3.17]{neukirchetal:2000}, corestriction is surjective by
  \cite[Proposition 3.3.8]{neukirchetal:2000}, and hence so is $f$. So
  all we need to do is verify that the formula we have given for $f$
  is equivariant with respect to $G$.  Note that $k_{\sigma\p} =
  k_{\p}-1$.  Hence, since $\inv_{\p}(\sigma c) =
  \inv_{\sigma^{-1}\p}(c)$, we have
        \begin{equation*}
                f(\sigma c) = f(c) - \sum_{\p\mid p}
                \inv_{\p}(c) \mathfrak{c} = f(c).
        \end{equation*}
\end{proof}

For a prime $\p$ of $L$ above $p$, denote by $(\ , \ )_{\p}$ the
Hilbert pairing on $L_{\p}^\times$ into $\mu_p$.  Let $\pi_L$
denote the projection map
$$
\pi_L \colon H^2(G_{K,S},\mu_p^{\otimes 2}) \to
\cl{K}/\Norm_{L/K}\cl{L} \otimes \mu_p. $$
\begin{theorem}
        \label{HCFpairformula}
        Let $\p_0$ be a prime of $L$ above $p$, and let $\mathfrak{c}$
        be as in Lemma~\ref{formulaforf}.  Let $a, b \in \McU{K}$ and
        suppose that $b = \Norm_{L/K}b'$ for some $b' \in \McU{L}$.
        Then
        \begin{equation}
        \label{eq:3}
                \pi_L \Mcpairing{a}{b} = \mathfrak{c}\tensor(a, N'b')_{\p_0},
                \quad N' = \sum_{k=1}^{p-1} k\sigma^k.
        \end{equation}
\end{theorem}

\begin{proof}
  As a standard property of the cup product, we have
  \begin{equation*}
    \Mcpairing{a}{b} = \cor_{L/K}(\res_{L/K}a,b')_{L,S}.
  \end{equation*}
  We evaluate $\pi_L\Mcpairing{a}{b}$ by taking the image of
  $(a,b')_{L,S}$ in $\bigoplus_{\p}H^2(L_{\p},
  \mu_p^{\otimes 2})$, namely $\bigoplus_{\p} (a,b')_{\p}$, and
  applying the map $f \otimes \mu_p$ to it.  By
  Lemma~\ref{formulaforf}, the result of this is
  \begin{equation}
        \pi_L\Mcpairing{a}{b} =
        \mathfrak{c} \otimes \sum_{\p\divides p}  (a,b')_{\p}^{k_{\p}}.
        \label{formula}
  \end{equation}
  Now, $(\sigma x,\sigma y)_{\sigma\p} = (x,y)_{\p}$ for any $x,
  y \in L^{\times}$. Thus~\eqref{eq:3} follows
  from~\eqref{formula} and the facts that $\sigma a = a$ and
  $\sigma^{-k_{\p}}\p_0 = \p$.
\end{proof}

We make the following observation on the applicability of
Theorem~\ref{HCFpairformula}.

\begin{lemma}
  If $|\McA{K}| = p$, then $\Norm_{L/K} \McU{L} = \McU{K}$.
\end{lemma}

\begin{proof}
  By Proposition~\ref{pairingandcapitulation}, we must show that
  $\McA{L}^G/\phi(\McA{K}) = 0$ (where $\phi$ is the natural map).
  Since $\McA{L}$ is finite and $G$ is cyclic, we have
  $$
  |\McA{L}^G| = |\McA{L}/(\sigma-1)|,
  $$
  and it follows from Lemma~\ref{genustheory} that this latter
  group has order $|\McA{K}|/p = 1$.
\end{proof}

We now focus on our main interest: $K = \QQ(\mu_p)$ and $S =
\{(1-\zeta)\}$, with $p$ an irregular prime satisfying Vandiver's
conjecture.  We consider the pairing $\blankeigpair{r}$, where $p$
divides $B_r$.  Now let $\alpha_{p-r}^p = \eta_{p-r}$, and set $L
= K(\alpha_{p-r})$, so that $\Delta$ acts on $G$ with eigenvalue
$\omega^{p-r}$.  Then $\pi_L$ amounts to the application of the
idempotent $\epsilon_{2-r}$.

We can exploit the action of $\Delta$ to simplify the computation
of $(N'b')_{\p_0}$ as follows.  Let $\iota_{\p}:L \into K_p =
\QQ_p(\mu_p)$ be the embedding corresponding to $\p$.  Let
$\Delta_0$ be the inertia group of $\p_0$, which we identify with
$\Gal(K_p/\QQ_p)$ via the projection onto $\Gal(K/\QQ)$.

\begin{proposition} \label{nontrivcondition}
  Let $3 \le i \le p-2$ be odd such that $p$ does not divide
  $B_{p-i}$.  Assume that $b = \eta_{r-i}$ is the norm of an $S$-unit
  $b' \in \McU{L}$, and choose $b'$ to have image in the
  $\omega^{p-r+i}$-eigenspace of $\McU{L}/\McO{L}^{\times p}$.  Then
  $\eigpair{\eta_i}{\eta_{r-i}}{r} \neq 0$
  if and only if $\iota_{\p_0} (N'b') \notin K_p^{\times p}$.
\end{proposition}

\begin{proof}
  Let $\delta \in \Delta_0$.  We see that
  \begin{equation*}
    \delta(N'b') = \prod_{k=1}^{p-1} \delta \sigma^k (b')^k =
    \prod_{k=1}^{p-1} \sigma^{k\omega^{p-r}(\delta)} \delta (b')^k.
  \end{equation*}
  Modulo $\McO{L}^{\times p}$, this is congruent to
  $$
    \prod_{k=1}^{p-1} \sigma^{k\omega^{p-r}(\delta)}(b')^{k\omega^{p-r+i}(\delta)}
    \equiv \prod_{k=1}^{p-1} \sigma^k (b')^{k\omega^{i}(\delta)}
    \equiv
    (N'b')^{\omega^{i}(\delta)}.
  $$
  Hence $N'b'$ has image in the $\omega^i$-eigenspace of
  $\McU{L}/\McO{L}^{\times p}$ under $\Delta_0$.

  By Theorem~\ref{HCFpairformula}, $\eigpair{\eta_i}{\eta_{r-i}}{r} = 0$
  if and only if $(\eta_i,N'b')_{\p_0} = 1$.  Since $p$
  does not divide $B_{p-i}$, the element $\eta_i$ is not locally a
  $p$th power.  Furthermore, we have seen that the elements $\eta_i$
  and $N'b'$ have image in the $\omega^{1-i}$ and $\omega^i$
  eigenspaces of $K_p^{\times}/K_p^{\times p}$, respectively.  These
  eigenspaces have dimension $1$ since $i \not\equiv 0,1 \pmod{p-1}$,
  and hence the result follows from the non-degeneracy and Galois
  equivariance of the norm residue symbol.
\end{proof}

For $p=37$, the condition of Proposition~\ref{nontrivcondition} is
computationally verifiable.

\begin{theorem}\label{NontrivialPairing}
  The pairing $\blankeigpair{32}$ for $p=37$ is nontrivial. Thus,
  Conjecture~\ref{k2conjecture} is true for $p=37$.
\end{theorem}

\begin{proof}[Sketch of Proof]
  Consider the fixed field $F$ of $\Delta_0$.  This is generated by
  the trace $x$ of a $p$th root $\alpha_{p-r}$ of (a choice of)
  $\eta_{p-r}$.  With the help of William Stein, we determined a
  minimal polynomial for $x$ by considering small primes $l$ that are
  primitive roots modulo $p$, computing the minimal polynomial of the
  image of $x$ in $\mathbb{F}_{l^{p-1}}[X]/(X^p-\eta_{p-r})$, and
  using the Chinese Remainder Theorem to find a $\QQ$-polynomial that
  $x$ satisfies.  Given this polynomial, Claus Fieker used Magma
  routines to compute the maximal order of $F$ and then a polynomial
  for $F$ with smaller discriminant (by far the most time-intensive
  steps), which made it possible to compute the $p$-unit group of $F$
  (by first computing the class group to ``sufficient precision").
  Now $F$ has two prime ideals above $p$, and the prime of $F$ below
  $\mathfrak{p}_0$ embeds $F$ into $\mathbb{Q}_p$.  We chose a
  $p$-unit that generates this prime, and this provided an element
  $b'$ as in Proposition~\ref{nontrivcondition} with $p =
  \Norm_{L/K}b'$.  We then computed the embeddings of $x$ at the
  primes $\mathfrak{p}_k$ from the embeddings of $\alpha_{p-r}$, which
  we obtained by factoring $X^p-\eta_{p-r}$ over $\QQ_p(\zeta)$.
  Writing $b'$ as a $\QQ$-polynomial in $x$, we then computed the
  image
  $$
  \iota_{\mathfrak p} (N'b') = \prod_{k=1}^{p-1} \iota_{\p_k}(b')^k
  $$
  to verify the condition of Proposition~\ref{nontrivcondition}.
  The Magma code is currently available at
  {\tt www.math.harvard.edu/\urltilde{}sharifi} and {\tt
  www.math.arizona.edu/\urltilde{}wmc}.
\end{proof}

Using the results of Section~\ref{idealclass}, we obtain the
following corollary, which implies, for example, that
$\mathbb{Q}(\!\! \sqrt[37]{37})$ has class number prime to $37$,
answering a question of Ralph Greenberg's.

\begin{corollary}
  Let $p = 37$, and let $L/K$ be a cyclic extension of degree $37$
  that is unramified outside $37$.  Then $|\pcl{L}| = 37$ if and only
  if $L$ is not contained in
  $\QQ(\zeta_{37^2},\alpha_{5},\alpha_{27})$, where $\alpha_{i}^{37} =
  \eta_{i}$ for any odd $i$.
\end{corollary}

\begin{proof}
  The values of the pairing $\blankeigpair{32}$ tell us that the
  subgroup of $D_K$ consisting of elements that pair trivially with
  all $37$-units of $K$ is $Q = \langle \zeta, \eta_5, \eta_{27}
  \rangle \cdot D_K^{37}$.  By Corollary~\ref{simplecase}, $|\pcl{L}|
  = 37$ if and only if $L = K(\alpha)$ with $\alpha^{37} \notin Q$.  The
  result now follows via Kummer theory.
\end{proof}

\section{Relations in the Galois group} \label{GaloisGroup}

Let us return to the general situation and notation of the
introduction with $n=p$, considering a free presentation
\eqref{presentation} of $\McG = G_{K,S}^{(p)}$.  The image of an
arbitrary relation $\Mcrho \in R$ in $\gr^2(\McF)/p \gr^1(\McF)$
was given by \eqref{rho}.  In $\gr^2 (\McF)$ itself, the relation
must have the form
\begin{equation} \label{fullrelation}
  \sum_{x \in X} a_x^{\Mcrho} p x + \sum_{x < x' \in X}
a_{x,x'}^{\Mcrho} [x,x'].
\end{equation}
To describe $a_x^{\Mcrho} \in \ZZ/p\ZZ$, we use the Bockstein
homomorphism
$$
B \colon H^1(\McG,\ZZ/p\ZZ) \rightarrow H^2(\McG,\ZZ/p\ZZ), $$
which is the coboundary in the long exact cohomology sequence of
$$
0 \to \ZZ/p\ZZ \xrightarrow{p} \ZZ/p^2\ZZ \to \ZZ/p\ZZ \to 0. $$
For $x \in X$ and $\Mcrho \in R$ we have
\cite[Proposition~3.9.14]{neukirchetal:2000}
\begin{equation} \label{BocksteinRelation}
        a_x^{\Mcrho} = -\rho(B(x^*)),
\end{equation}
where $x^* \in X^*$ is the dual to $x$.

\begin{lemma} \label{Bockstein}
  Let $K$ be a number field containing $\mu_p$ and $S$ a set of primes
  containing those above $p$.  For $a \in D_K$ the homomorphism $\Phi
  = B \otimes \mathrm{id}_{\mu_p}$ is given, abusing notation, by
  $$
  \Phi(a) \otimes \zeta = \mathfrak{a} \otimes \zeta -
  (\zeta,a)_S,$$
  where $\mathfrak{a}^p = a\McO{K}$.
\end{lemma}

\begin{proof}
  By comparison with the Kummer sequence \eqref{kummer} for
  $n = p$, the coboundary map $B^*$ in the cohomology of the short
  exact sequence
  $$1 \to \mu_p \to \mu_{p^2} \xrightarrow{p} \mu_p \to 1,$$
  is seen
  to be given by $B^*(a) = \mathfrak{a} \pmod{p}$.  We compute $B -
  (B^* \otimes j)$, where $j$ is the identity map on
  $\mu_p^{\otimes(-1)}$, in terms of cocycles.

  On the one hand, for $f \in H^1(\McG,\ZZ/p\ZZ)$, we have that $B(f)$
  is the class of
  $$
  (\sigma,\tau) \mapsto \frac{1}{p}(\tilde{f}(\tau) +
  \tilde{f}(\sigma) - \tilde{f}(\sigma\tau))
  $$
  for an arbitrary lift of $f$ to a map $\tilde{f} \colon \McG \to
  \ZZ/p^2\ZZ$.  On the other hand, $B^* \otimes j$ takes $f$ to the
  class of
  $$
  (\sigma,\tau) \mapsto \frac{1}{p}(\chi(\sigma)\tilde{f}(\tau) +
  \tilde{f}(\sigma) - \tilde{f}(\sigma\tau)),
  $$
  where $\chi \colon \McG \to (\ZZ/p^2\ZZ)^{\times}$ is the
  cyclotomic character associated with a root of unity of order $p^2$.
  The difference of these cocycles is
  $$
  (\sigma,\tau) \mapsto -\frac{\chi(\sigma)-1}{p} f(\tau),
  $$
  and here $\sigma \mapsto \zeta^{(\chi(\sigma)-1)/p}$ is the
  Kummer character associated with $\zeta$.  By the well-known formula
  for the cup product of two homomorphisms as their product, we have
  the result.
\end{proof}

Again, let us focus on $K = \QQ(\zeta_p)$ and $S = \{ (1-\zeta)
\}$, with $p$ satisfying Vandiver's conjecture.  We describe a
minimal generating set $X$ of $\McG$.  Let $M$ denote the set of
integers $m$ with $2 \le m \le p$ and either $m$ odd or
$(A_K/p)(\omega^{p-m})$ nontrivial.  (We take the given interval,
instead of $1 \le m \le p-1$, for compatibility with
Section~\ref{fundamental}.)  We let $M_{\mathrm{o}}$ and
$M_{\mathrm{e}}$ denote the odd and even elements of $M$,
respectively.  For each $m \in M \cup \{0\}$, we choose an element
$x_m \in \McG$ with image generating the $\omega^m$-eigenspace of
$\gr^1(\McG)$, subject to the following normalizations.  For $m
\in M_{\mathrm{o}}$, we assume that $x_m \alpha_m = \zeta
\alpha_m$ for $\alpha_m$ a $p$th root of $\eta_m$, the cyclotomic
unit defined in Section~\ref{cyclpunits}.  For $m \in
M_{\mathrm{e}}$, let $\mathfrak{b}_{m} \in A_K(\omega^{p-m})$ be
such that $\mathfrak{b}_{m} \tensor \zeta$ maps to $1$ under the
isomorphism \eqref{isochoice} (chosen in defining the pairing
$\blankeigpair{m}$), and let $f_{p-m}$ be the Iwasawa power series
for $A_{\infty}(\omega^{p-m})$ with
$$
f_{p-m}((1+p)^s-1) = L_p(\omega^m,s), $$ for $s \in \ZZ_p$, where
$L_p(\omega^m,s)$ is the $p$-adic $L$-function.  Choose $b_m \in
D_K$ with image in $(D_K/D_K^p)(\omega^{p-m})$ such that
\begin{equation} \label{normalize}
        b_m \mapsto \mathfrak{b}_m^{f_{p-m}(0)/p} \qquad \text{under\ } D_K \to A_K[p].
\end{equation}
Writing $b_m = \beta_m^p$, we require that $x_m(\beta_m) =
\zeta\beta_m$.  Finally, for $m = 0$, we let $x_0 = \gamma$
satisfy $\gamma(\xi) = \xi^{1+p}$ for all $\xi \in
\mu_{p^{\infty}}$.

Let
$$
X = \{ x_m \colon m \in M \cup \{0\} \}.
$$
Then $X$ provides a dual basis to the basis of $H^1(\McG,\mu_p)$
given by the elements $\eta_m$ for $m \in M_{\mathrm{o}}$, $b_m$
for $m \in M_{\mathrm{e}}$, and $\zeta$ (under the isomorphism
$\mu_p \xrightarrow{\sim} \ZZ/p\ZZ$ provided by $\zeta$).

Now fix $r \in M_{\mathrm{e}}$.  Recall that $e_{i,r}$ was defined
to be $\eigpair{\eta_i}{\eta_{r-i}}{r}$ for $i \in
M_{\mathrm{o}}$.  If $i \in M_{\mathrm{e}}$ and the least positive
residue $j$ of $r-i$ modulo $p-1$ is also in $M_{\mathrm{e}}$ (so
that $p$ divides $B_i$, $B_{j}$ and $B_r$!) we then set $e_{i,r} =
\eigpair{b_i}{b_j}{r}$.

Identify $1+T$ with the restriction of $\gamma$ to $K_{\infty} =
\QQ(\mu_{p^{\infty}})$.  Define $g_{p-r} \in \Lambda \iso
\ZZ_p[[T]]$ by the relation
\begin{equation} \label{duality}
        g_{p-r}((1+p)^s-1) = f_{p-r}((1+p)^{1-s}-1)
\end{equation}
for every $s \in \zp$.

\begin{theorem} \label{grouprelation}
  For $r \in M_{\mathrm{e}}$ with $p$ satisfying Vandiver's
  conjecture, there is a relation in $\gr^2(\McG)$ of the form
  $$
  g_{p-r}(0) x_r + g_{p-r}'(0) [\gamma,x_r] + \sum_{\substack{i<j
      \in M \\ i+j \equiv r \pmod{p-1}}} e_{i,r} [x_i,x_j] = 0.
  $$
\end{theorem}

\begin{proof}
  Choose a relation $\rho \in \McF$ such that the map
  $H^2(\McG,\mathbb{Z}/p\mathbb{Z}) \to \mathbb{Z}/p\mathbb{Z}$
  corresponding to $\rho$ factors through the negative of the
  isomorphism \eqref{isochoice} used in defining $\blankeigpair{r}$.
  By the expression \eqref{fullrelation} for $\rho$, the identities
  \eqref{BasicRelation} and \eqref{BocksteinRelation}, and eigenspace
  considerations, we see that
  $$
  \rho \equiv a_{x_r}^{\rho} p x_r +
  a_{\gamma,x_r}^{\rho} [\gamma,x_r] + \sum_{\substack{i<j \in M \\
      i+j \equiv r \pmod{p-1}}} a_{x_i,x_j}^{\rho} [x_i,x_j]
  \pmod{\mathrm{Fil}^3 \McF},
  $$
  with $a_{x_i,x_j}^{\rho} = e_{i,r}$.

  Let $f = f_{p-r}$ and $g = g_{p-r}$.  We claim that
  \begin{eqnarray*}
     a_{x_r}^{\Mcrho} = g(0)/p \pmod{p} & \mathrm{and} &
     a_{\gamma,x_r}^{\Mcrho} = g'(0) \pmod{p}.
  \end{eqnarray*}
  Recalling the notation and statement of \eqref{normalize},
  Proposition~\ref{zetapairing} and Lemma~\ref{Bockstein} imply that
  \begin{eqnarray*}
     B(x_r^*) \otimes \zeta = \mathfrak{b}_r^{f(0)/p + f'(0)} &
     \mathrm{and} &
     (\zeta,b_r)_S = \mathfrak{b}_r^{-f'(0)} \otimes \zeta.
  \end{eqnarray*}
  Applying $-\rho$, we obtain by \eqref{BocksteinRelation} and
  \eqref{BasicRelation} that
  \begin{eqnarray*}
     a_{x_r}^{\Mcrho} = f(0)/{p} + f'(0) \pmod{p} & \mathrm{and} &
     a_{\gamma,x_r}^{\Mcrho} = -f'(0) \pmod{p}.
  \end{eqnarray*}
  That these agree with $g(0)/p$ and $g'(0)$ follows from the
  definition \eqref{duality} of $g$ in terms of $f$.
\end{proof}

From our table of the pairings and some basic Bernoulli number
computations, we obtain the following corollary of
Theorem~\ref{NontrivialPairing} (for a particular choice of
$x_{32}$).

\begin{corollary} \label{37relation}
  For $p = 37$ and $r = 32$, there is a relation in $\gr^2(\McG)$ of
  the form
        \begin{multline*}
          37 y - 3 [\gamma,y] - 11[x_3,x_{29}] - [x_7,x_{25}] +
          [x_9,x_{23}] - 2[x_{11},x_{21}] - 6[x_{13},x_{19}] \\
          - 3[x_{15},x_{17}] - [x_{31},x_{37}] + 11[x_{33},x_{35}] = 0
        \end{multline*}
        with $y = x_{32}^c$ for some $c \in (\ZZ/37\ZZ)^{\times}$.
\end{corollary}

\begin{proof}
  The coefficients of all but the first two terms are obtained from
  the table of calculated pairings.  It is easy to check that (see the
  proof of \cite[Corollary 10.17]{washington:1997})
  $$
  pf_{p-r}'(0) \equiv \frac{B_r}{r} - \frac{B_{r+p-1}}{r-1} \equiv
  16p \pmod{p^2} $$
  and
  $$
  f_{p-r}(0) = \frac{r-2}{r} B_r - B_{r+p-1} \equiv 14p \pmod{p^2},
  $$
  in order to compute the first two coefficients (up to a scalar
  relative to the others).
\end{proof}

\section{Relationship with pro-$p$ fundamental groups}
\label{fundamental}

We consider the curve $V = {\mathbb P}^1 - \{0,1,\infty\}$ over
$\QQ$, and let $\overline V = V \times_\QQ \overline \QQ$.  The
natural identification of the absolute Galois group $G_\QQ$ with
$\Aut(\overline V/V)$ induces a representation
$$
\phi \colon G_{\QQ} \to \mathrm{Out}(\pi_1^{(p)}(\overline V))
$$
where $\pi_1^{(p)}(\overline V)$ denotes the pro-$p$ fundamental
group.  As before, let $K = \QQ(\mu_p)$ and $S = \{ (1-\zeta) \}$.
Then $\phi$ factors through $\McG$, and one can put a ``weight''
filtration $F^m\McG$ on $\McG$, consisting of the subgroups
corresponding to the fixed fields of the various induced
representations
$$
\phi_m \colon \McG \to \mathrm{Out}(\pi_1/\pi_1(m+1)),
$$
where $\pi_1 = \pi_1^{(p)}(\overline V)$ and $\pi_1(m+1)$ is the
$(m+1)$th term in the descending
central series of $\pi_1$. In particular, the fixed field of
$\phi_1$ is $K_{\infty} = \QQ(\mu_{p^{\infty}})$. Consider the
graded $\zp$-Lie algebra
$$
\g = \oplus_{m=1}^{\infty} \gr^m \g,\quad\mathrm{where}\quad\gr^m
\mathfrak{g} = F^m\McG/F^{m+1}\McG.
$$
Then $\gr^m \g$ is torsion-free of finite $\ZZ_p$-rank, and
$G_{\QQ}$ acts on it by the $m$th power of the $p$-adic cyclotomic
character \cite{ihara:1986}.

Let $G$ denote the closed normal subgroup of $\McG$ with fixed
field $K_{\infty}$.  Note that $F^mG = F^m\McG$ for $m \ge 1$.
For odd $m \ge 1$, let $\kappa_m \colon G \to \ZZ_p(m)$ denote the
$\Gal(K_{\infty}/\QQ)$-equivariant homomorphism dual to an
appropriate sequence of cyclotomic $p$-units as defined in
\cite{ihara:1999} and known to be nontrivial by \cite{soule:1981}.
Then $\kappa_m$ induces a nontrivial map $\kappa_m \colon \gr^m \g
\to \ZZ_p$ for odd $m \ge 3$ \cite[Proposition 1]{ihara:1989}.
Let $\tilde{\sigma}_m$ denote an element of $F^m G$ such that
$v_p(\kappa_m(\tilde{\sigma}_m))$ is minimal.  Then
$\tilde{\sigma}_m$ restricts to a nontrivial element $\sigma_m \in
\gr^m \g$.

Let $\mathfrak{h}$ denote the Lie subalgebra of $\mathfrak{g}$
generated by the $\sigma_m$.  Hain and Matsumoto have proven a
conjecture of Deligne that $\g \subset \mathfrak{h} \otimes \QQ_p$
\cite{hain-matsumoto:2001} (see also \cite[Section
3.7]{goncharov:2001} for a description of motivic arguments of
Beilinson and Deligne that lead to this result). Deligne has
further conjectured that $\mathfrak{g} \otimes \QQ_p$ is free on
the $\sigma_m$.  On the other hand, Sharifi \cite[Theorem
1.3]{sharifi:2001} has shown that Greenberg's conjecture, as
described in Section~\ref{greenberg}, implies that $\mathfrak{g}$
itself is not free on the $\sigma_m$ if $p$ is irregular.
Deligne's conjecture would then imply that $\mathfrak{h} \neq
\mathfrak{g}$. In this section, we describe relations in
$\mathfrak{g}$, conjecturally nontrivial, in terms of the cup
product $\blankpairing$ (Theorem~\ref{lierelation}). Taking $p =
691$ as an example, we
see in Theorem~\ref{p691m12} that the nontriviality of
$\blankeigpair{12}$ is, in fact, equivalent to $\gr^{12}
\mathfrak{h} \neq \gr^{12} \mathfrak{g}$.

Let $\mathfrak{s}$ the free $\ZZ_p$-Lie algebra on generators
$s_i$ with $i \ge 3$ odd, and let $\psi \colon \mathfrak{s} \to
\g$ denote the map given by $s_i \mapsto \sigma_i$ for each $i$.
Deligne's conjecture is equivalent to the statement that $\psi$ is
injective. We say that Deligne's conjecture holds in degree $i$ if
$\gr^i \psi$ is injective.

Assume $p$ satisfies Vandiver's Conjecture for the remainder of
the section. Recall the notation of Section~\ref{GaloisGroup}. Let
$\Delta$ be a choice of lift of $\Gal(K/\QQ)$ to a subgroup of
$\McG$ of order $p-1$.  We may choose $x_m$ for each $m \in M \cup
\{0\}$ such that
$$
\delta x_m \delta^{-1} = x_m^{\omega(\delta)^m} $$ for each
$\delta \in \Delta$ by \cite[Lemma 2.1]{sharifi:2001}, which
immediately forces $x_m \in F^m G$ for $m \in M$.  Furthermore,
for $m \in M_{\mathrm{o}}$, we may take $\tilde{\sigma}_m = x_m$.
For each $m \in M$, we let $\bar{x}_m$ denote the image of $x_m$
in $\gr^m \g$, which for even $m$ may or may not be trivial.

Recall that $\Fil^k \McG$ was defined as the descending
$p$-central series of $\McG$. We define the induced filtration on
$\g$:
$$
  {\underline \Fil}^k \g = \bigoplus_{m=1}^\infty \frac{\Fil^k \McG
  \cap F^m G}{\Fil^k \McG \cap F^{m+1} G},
$$
which differs from the descending central $p$-series on $\g$.

\begin{theorem} \label{lierelation}
  Let $m \in M_{\mathrm{e}}$ for $p$ satisfying Vandiver's conjecture.
  Then there is a relation in $\gr^m \g$ of the form
  \begin{equation} \label{therelation}
    B_m \bar{x}_m \equiv m \sum_{\substack{i<j \in M \\ i+j = m}}
    e_{i,m}[\bar{x}_i,\bar{x}_j] \pmod{\gr^m{\underline \Fil}^3\g},
  \end{equation}
  where the $e_{i,m}$ are as defined in Section~\ref{GaloisGroup}.
\end{theorem}

\begin{proof}
  We note that for $x \in F^i G$, we have
  \begin{equation} \label{actbygamma}
        x^{\gamma-1} \equiv x^{(1+p)^i -1} \bmod F^{i+1}G,
  \end{equation}
  since $\gr^i \g$ has Tate twist $i$.  Let $g = g_{p-m}$ be as in
  \eqref{duality}.  Since
  $$
  g((1+p)^m-1) = L_p(\omega^m,1-m) = (p^m-1)\frac{B_m}{m},
  $$
  applying \eqref{actbygamma}, we obtain
  \begin{equation*}
    x_m^{g(T)} \equiv x_m^{(p^m-1)B_m/m} \pmod{F^{m+1}G}.
  \end{equation*}
  The result now follows from Theorem~\ref{grouprelation} by reducing
  its relation modulo the image of $F^{m+1} G$ in $\gr^2(\McG)$.
\end{proof}

We derive some consequences of this result regarding the freeness and
generation of $\g$.  In order to do so, we must compare the filtration
${\underline \Fil}^{\cdot} \g$ on $\g$ with the descending central
$p$-series on the simpler Lie algebra $\h$. This will proceed in
several steps. We begin with the following lemma.

\begin{lemma} \label{freesubgroup}
  Let $\mathcal{H}$ be the pro-$p$ subgroup of $\McG$ generated by the
  $x_i$ for $i \in M_{\mathrm{o}} \cup \{ 0 \}$.  Let $H$ be the
  pro-$p$ subgroup of $G$ generated by the $\tilde{\sigma}_i$ for odd
  $i \ge 3$.  Then $\mathcal{H}$ and $H$ are freely generated as
  pro-$p$ groups on these sets of elements.
\end{lemma}

\begin{proof}
  Hain and Matsumoto \cite[Theorems 7.3, 7.4]{hain-matsumoto:2001}
  have demonstrated the existence of a filtration on $G$ with graded
  quotient a $\ZZ_p$-Lie algebra that injects into a free graded
  $\QQ_p$-Lie algebra on the images of the $\tilde{\sigma}_i$ in
  degree $i$.  As remarked by Ihara \cite[Section 6]{ihara:1999}, this
  implies that $H$ must be free on the $\tilde{\sigma}_i$.  By
  \cite[Lemma 3.1c]{sharifi:2001} (as in the proof of Theorem 1.3
  therein), this implies the desired freeness on the generators of
  $\mathcal{H}$.
\end{proof}

For $Z$ a pro-$p$ group or $\ZZ_p$-Lie algebra, we let $Z(k)$
denote the $k$th term in its descending central series and $\Fil^k
Z$ the $k$th term in its descending central $p$-series.
Next, we compare $\Fil^{\cdot} \h$ with the filtration induced on
$\h$ by $\Fil^{\cdot} H$.
Note that we use $H$, as opposed to
$\mathcal{H}$ (at this point), since the weight filtration is
defined by a filtration on $H$.
\begin{lemma} \label{lowercentral}
  Let $m$ and $k$ be positive integers.  If $\gr^m \psi$ is injective
  for $i < m$ then
  $$
  \gr^m \h(k) \iso \frac{H(k) \cap F^m G}{H(k) \cap F^{m+1} G}
  $$
  and
  $$
  \gr^m \Fil^k \h \iso \frac{\Fil^k H \cap F^m G}{\Fil^k H \cap
    F^{m+1} G}.
  $$
\end{lemma}

\begin{proof}
  By definition of $\h$, we may lift any element of $\gr^m \h(k)$ to
  an element of $H(k) \cap F^m G$.  We must show that, conversely,
  an element of
  $H(k)\cap F^m G$ projects to an element of $\gr^m \h(k)$.  By the
  injectivity of $\gr^i \psi$ in weights $i < m$, the group $H/(H \cap
  F^m G)$ is isomorphic to the free pro-$p$ subgroup on the
  $\tilde{\sigma}_i$ modulo the pro-$p$ subgroup generated by
  commutators
  $[\tilde{\sigma}_{m_1},\ldots[\tilde{\sigma}_{m_{j-1}},\tilde{\sigma}_{m_j}]\ldots]$
  with $\sum m_i \ge m$ \cite[Section 6]{ihara:1999}.  Using the
  freeness of $H$ in Lemma~\ref{freesubgroup}, $H(k) \cap F^m G$ is
  then the normal pro-$p$ subgroup of $H$ generated by those among the
  above commutators with $j=k$, which clearly project to elements of
  $\gr^m \h(k)$.  The same arguments hold with the descending
  central series terms replaced by descending central $p$-series terms.
\end{proof}

Now we compare the filtration ${\underline \Fil}^{\cdot} \g$
induced by $\Fil^{\cdot} \McG$ with the filtration on $\h$ induced
by $\Fil^{\cdot} \mathcal{H}$, since $\McG$ is more closely
related to $\mathcal{H}$ than to $H$.

\begin{lemma} \label{1stfiltrationlemma}
  Let $m$ and $k$ be positive integers with $m \le p+r-2$,
  with $r$ the minimal element of $M_{\mathrm{e}}$.
  If $\gr^i
  \psi$ is
surjective
  for $i \le m$, then
  $$
  \gr^m {\underline \Fil}^k \g
\iso
  \frac{\Fil^k \mathcal{H} \cap F^{m} G}{\Fil^k \mathcal{H} \cap F^{m+1} G}.
  $$
\end{lemma}

\begin{proof}
  We first show that $x_r \in [H,H] \cdot
  F^{m+1} G$ if $r \in M_{\mathrm{e}}$ with $r \le m$.  Note that
  since $r$ is even, $\gr^r \h = \gr^r [\h,\h]$ by definition.  Since
  $x_r \in F^r G$ and $\gr^r \psi$ is surjective, we have $\bar{x}_r
  \in \gr^r [\h,\h]$.
  Thus, we obtain
  $x_r \in [H,H] \cdot F^{r+1} G$.
  Now assume that $x_r \in [H,H] \cdot
  F^l G$ for some $l \not\equiv r \pmod{p-1}$ with $l \le m$.  We
  remark that
  $$
  \frac{[H,H] \cdot F^l G}{[H,H] \cdot F^{l+1} G} \iso \frac{F^l
    G}{([H,H] \cap F^l G) \cdot F^{l+1} G}
    \iso \gr^l \h^{\ab},
  $$
  and $\Gal(K_{\infty}/\QQ)$ acts on the latter group by the $l$th
  power of the cyclotomic character.  On the other hand,
  $$\delta x_r \delta^{-1} = x_r^{\omega(\delta)^r} \qquad \text{\
    for\ } \delta \in \Delta,$$
  and this forces $x_r \in [H,H] \cdot
  F^{l+1} G$.  By recursion, since $m < r + p - 1$, we have $x_r \in
  [H,H] \cdot F^{m+1} G$.

  For a tuple ${\bf m} = (m_1, \ldots, m_j)$ with $m_i \in M \cup
  \{0\}$ for each $i$, and $m_{j-1} < m_{j}$ if $j \ge 2$, let $$
  x_{{\bf m}} = [x_{m_1},\ldots,[x_{m_{j-1}},x_{m_{j}}]\ldots]. $$
  Then
  $\Fil^k \McG$ is generated as a normal subgroup by those elements of
  the form $p^{k-j} x_{{\bf m}}$, with ${\bf m}$ a tuple of length $j
  \le k$.  Since $x_r \in [H,H] \cdot F^{m+1} G$ for $r \in
  M_{\mathrm{e}}$, the images of those elements $p^{k-j}x_{{\bf m}}$
  with some $m_i \in M_{\mathrm{e}}$ are redundant as elements of the
  induced generating set of $\Fil^k \McG/(\Fil^k \McG \cap F^{m+1}
  G)$, and hence $\Fil^k \mathcal{H}$ surjects onto the latter
  quotient, finishing the proof.
\end{proof}

In the following proposition, we conclude our discussion of
filtrations by filling in the intermediate comparison between the
filtrations on $\h$ induced by $\Fil^{\cdot} H$ and $\Fil^{\cdot}
\mathcal{H}$. Note that these filtrations will always disagree in
sufficiently large weight, since for each $k$, there exists $i$
sufficiently large such that $\tilde{\sigma}_i \in \Fil^k
\mathcal{H}$.
\begin{proposition} \label{filtrationlemma}
  Let $m$ and $k$ be positive integers with $m \le p+1$.
  If $\gr^i
  \psi$ is bijective for $i \le m$, then
  $$
  \gr^m {\underline \Fil}^k \g = \gr^m \Fil^k \h.
  $$
\end{proposition}

\begin{proof}
  By the second isomorphism in Lemma~\ref{lowercentral} and the
  isomorphism of Lemma~\ref{1stfiltrationlemma}, it suffices to show
  that
  $$
  \frac{\Fil^k H \cap F^m G}{\Fil^k H \cap F^{m+1} G} \iso
  \frac{\Fil^k \mathcal{H} \cap F^m G}{\Fil^k \mathcal{H} \cap F^{m+1}
    G}.
  $$
  We clearly have that $\Fil^k H \subseteq \Fil^k \mathcal{H}$, and
  we are left to verify that
  \begin{equation} \label{pseriesclaim}
        (\Fil^k \mathcal{H} \cap F^{m} G) \cdot F^{m+1} G \subseteq
        (\Fil^k H \cap F^{m} G) \cdot F^{m+1} G
  \end{equation}

  From the generating set $\{x_i \colon i \in M_{\mathrm{o}} \cup
  \{0\}\}$ of $\mathcal{H}$, we may define a generating set $\{
  x_{i,n} \colon i \in M_{\mathrm{o}},\ n \ge 0 \}$ of $H$, taking
  $x_{i,0} = x_i$ and
  $$
    x_{i,n+1} = \gamma x_{i,n} \gamma^{-1} x_{i,n}^{-(1+p)^{i+n(p-1)}}.
  $$
  By \cite[Lemma 3.1b]{sharifi:2001}, the $x_{i,n}$ freely generate $H$
  as a pro-$p$ group.
  Note that $x_{i,n} \in \Fil^{n+1} \mathcal{H} -
  \Fil^{n+2} \mathcal{H}$ by Lemma~\ref{freesubgroup}.  Thus $H \cap
  \Fil^k \mathcal{H}$ is generated as a normal pro-$p$ subgroup by
  those elements of the form
  $$
  p^{k-J}[x_{m_1,n_1},\ldots,[x_{m_{j-1},n_{j-1}},x_{m_j,n_j}]\ldots]
  $$
  with $m_t \in M_{\mathrm{o}}$, $n_t \ge 0$, and
  $J = j+ \sum n_t$
  for $1 \le t \le j$.  Furthermore, it follows as in
\cite[Lemma
  2.2]{sharifi:2001}, that the elements $\tilde{\sigma}_i$ with $i \ge
  3$ odd may be chosen such that
  $$
  x_{i,n}\tilde{\sigma}_{i+n(p-1)}^{-1} \in [H,H] \cap \Fil^{n+1}
  \mathcal{H}. $$
  Therefore, $H
  \cap \Fil^k \mathcal{H}$ is generated as a normal pro-$p$ subgroup
  by the elements
  $$
  p^{k-J}[\tilde{\sigma}_{m_1+n_1(p-1)},\ldots,
  [\tilde{\sigma}_{m_{j-1}+n_{j-1}(p-1)},\tilde{\sigma}_{m_j+n_j(p-1)}]\ldots],
  $$
  with $m_t$, $n_t$ and $J$ as before.  Since
  $\tilde{\sigma}_{i+n(p-1)} \in F^{m+1}G$ if $n \ge 1$, to show
  \eqref{pseriesclaim}, it suffices to verify that
  $$
  p^{k-j}[\tilde{\sigma}_{m_1},\ldots,
  [\tilde{\sigma}_{m_{j-1}},\tilde{\sigma}_{m_j}]\ldots] \in \Fil^k H,
  $$
  but this is true by definition.
\end{proof}

Combining Proposition \ref{filtrationlemma} with
Theorem~\ref{lierelation}, we obtain the following.

\begin{proposition} \label{notfree}
  Assume that $e_{i,m}$ is nonzero for some $i \in M_{\mathrm{o}}$
  with $i < m/2$ and $m \in M_{\mathrm{e}}$.  Then $\g$ is not freely
  generated by the elements $\sigma_i$.  In fact, there is an $i \le
  m$ for which $\gr^i \psi$ is not is an isomorphism.
\end{proposition}

\begin{proof}
  Assume that $\gr^i \psi$ is bijective for all $i \le m$, so the
  relation \eqref{therelation} holds modulo $\gr^m \Fil^3 \h$ by
  Proposition \ref{filtrationlemma}.  Then $B_m \bar{x}_m \in p \cdot
  \gr^m \h$ by the surjectivity of $\gr^m \psi$.  Furthermore, for any
  $r \in M_{\mathrm{e}}$ with $r < m$, we have $\bar{x}_r \in [\h,\h]$
  by the surjectivity of $\gr^r \psi$, as in the proof of Lemma
  \ref{1stfiltrationlemma}.  Since some $e_{i,m} \neq 0$, reducing
  \eqref{therelation} modulo $p\h + \h(3)$ exhibits a contradiction of
  the injectivity of $\gr^m \psi$.
\end{proof}

This can be improved as follows, when $m$ is the smallest positive
even integer such that $p$ divides $B_m$.

\begin{theorem} \label{notgenerated}
  Let $m$ be the minimal element of $M_{\mathrm{e}}$ for an irregular
  prime $p$ satisfying Vandiver's conjecture.  Assume that Deligne's
  conjecture holds in degrees $i \le m$.  Then $\gr^m \g =
  \ZZ_p\bar{x}_m + \gr^m \h$ and $\gr^m \g^{\ab}$ is generated by the
  image of $\bar{x}_m$.  Furthermore, if $e_{i,m}$ is nonzero for some
  $i \in M_{\mathrm{o}}$ with $i < m/2$, then in fact $\gr^m \h
  \subsetneq \gr^m \g$ and $\gr^m \g^{\ab}$ is nontrivial.
\end{theorem}

\begin{proof}
  First, we remark that $\gr^i \psi$ is
  not only injective, but bijective in degrees $i < m$
  by \cite[Theorem I.2(ii)]{ihara:1999}.  From this, it is easy to see
  that $\bar{x}_m$ and $\gr^m \h$ generate $\gr^m \g$.  That is,
  $\gr^m [\g,\g] = \gr^m [\h,\h]$ by the bijectivity in lower degrees,
  and $\gr^m \g^{\mathrm{ab}}$ is generated by the image of
  $\bar{x}_m$ as shown, e.g, in the proof of
  \cite[Theorem 4.1]{sharifi:2001}.
  By Proposition \ref{notfree}, we know that $\gr^m \psi$ is not
  bijective, hence not surjective, finishing the proof.
\end{proof}

Let $\McD$ denote Ihara's stable derivation algebra
\cite{ihara:1999}, which is a graded Lie algebra over $\ZZ$.  More
specifically, it is a Lie algebra of derivations of the free
graded Lie algebra $\mathcal{F}$ on two variables $x$ and $y$ over
$\ZZ$ and consists of $D \in \gr^m \McD$ such that $D(x) = 0$ and
$D(y) = [y,f_D]$ with $f_D \in \gr^m \mathcal{F}$ satisfying
certain relations.

Ihara has shown that there is an injection of graded $\ZZ_p$-Lie
algebras \cite{ihara:1989}
\[ \iota \colon \g \to \McD \otimes \ZZ_p. \]
He has also made the following conjecture.
\begin{conjecture}[Ihara] \label{iotaconj}
  The map $\gr^m \iota$ is an isomorphism for $m < p$.
\end{conjecture}

Consider the map
\[ \lambda_m \colon \gr^m \McD \to \ZZ \]
(denoted $\gr^m(c)$ in \cite{ihara:1999}) given by
$$
f_D \equiv \lambda_m(D)[x,[x,\ldots[x,y]\ldots]]
\pmod{\mathrm{terms\ of\ degree\ } \ge 2\mathrm{\ in\ }y}.
$$
We remark that $\lambda_m = 0$ if $m$ is even or $m = 1$.
Extending $\lambda_m$ $\zp$-linearly, we may precompose with
$\gr^m\iota$ to obtain a map $\lambda_m^{(p)}$ satisfying the
formula \cite{ihara:1989,ihara:1999}
\begin{equation} \label{relate}
    \kappa_m = (p^{m-1}-1)(m-1)!\lambda_m^{(p)}
\end{equation}
for odd $m \ge 3$.

\begin{lemma} \label{Dmcondition}
  For $m \in M_{\mathrm{o}}$, there exists $D_m \in \gr^m \McD$
  such that $\lambda_m(D_m)$ is the positive generator of the image of
  $\lambda_m$ and
  \begin{equation} \label{D_m}
    D_m \equiv -(m-1)!\lambda_m(D_m)\iota(\sigma_m) \pmod{p\McD
    \otimes \zp}.
  \end{equation}
\end{lemma}

\begin{proof}
  We need only show that
the two defining properties of $D_m$ are consistent.
  If $m \le p$, we have
  \[
    \lambda_m^{(p)}(\sigma_m) \equiv -1/(m-1)! \pmod{p}
  \]
  by equation \eqref{relate}, and consistency follows from
  applying $\lambda_m$ to both sides of \eqref{D_m}.
\end{proof}

Ihara has conjectured the existence of $p$-torsion in the $m$th
graded piece of $\McD^{\ab}$ when $p$ divides $B_m$
\cite[Conjecture II.2]{ihara:1999}.  We will focus on a case in
which the stable derivation algebra has been calculated
sufficiently to allow comparison with the relation in
Theorem~\ref{lierelation}.  That is, when $p = 691$ and $m = 12$,
Ihara has exhibited a relation
\begin{equation} \label{derivation12}
  691\delta = 2[D_3,D_9] - 27[D_5,D_7]
\end{equation}
for some $\delta \in \gr^{12} \McD$.  On the other hand, Matsumoto
has shown that $[D_3,D_9]$ and $[D_5,D_7]$ form a basis of
$\gr^{12} \McD \otimes \QQ$ and generate $\gr^{12} [\McD,\McD]$
\cite[Appendix A]{matsumoto:1995}, which implies that $\delta
\notin [\McD,\McD]$. Furthermore, he has verified that the image
of $\delta$ generates $\gr^{12} \McD^{\ab}$.

\begin{proposition} \label{equivalence}
  For $p = 691$, Conjecture~\ref{iotaconj} in degree $m = 12$ is
  equivalent to the statement that $\gr^{12} \h \subsetneq \gr^{12}
  \g$.
\end{proposition}

\begin{proof}
  We remark that $\gr^i \psi$ is injective for $i \le 11$, as follows
  directly from the main results of \cite{ihara:1989}.  By
  Theorem~\ref{notgenerated}, the top row of the commutative diagram
  \begin{equation*}
        \xymatrix{0 \ar[r] & \gr^{12} [\h,\h] \ar[r]\ar[d] & \gr^{12} \g
        \ar[r]\ar[d]^{\gr^{12} \iota} & \gr^{12} \g^{\ab} \ar[r]\ar[d] & 0 \\
        0 \ar[r] & \gr^{12} [\McD,\McD] \otimes \ZZ_{691} \ar[r] & \gr^{12}
        \McD \otimes \ZZ_{691} \ar[r] & \gr^{12} (\McD \otimes \ZZ_{691})^{\ab} \ar[r] & 0 \\
        }
  \end{equation*}
  is exact and $\gr^{12} \g^{\ab}$ is nonzero if and only if
  $\gr^{12} \h \subsetneq \gr^{12} \g$.
The bottom row is exact by definition.
  Since $\gr^{12}[\McD,\McD]$ is generated by $[D_3,D_9]$ and $[D_5,D_7]$,
  Lemma~\ref{Dmcondition} implies that the leftmost vertical arrow is
  a surjection.  Since $\iota$ is injective, the rightmost vertical
  arrow is now forced to be an injection as well.  Furthermore, noting
  \eqref{derivation12} and
Matsumoto's results discussed after it, we have
  that $\gr^{12} \McD^{\ab}$ is cyclic
  of order $691$.  Therefore, $\gr^{12} \iota$ is an isomorphism if
  and only if $\gr^{12} \g^{\ab}$ is nonzero, hence the result.
\end{proof}

By Proposition~\ref{equivalence}, the following shows that
Conjectures~\ref{k2conjecture} and~\ref{iotaconj} are equivalent
for the irregular pair $(691,12)$.

\begin{theorem} \label{p691m12}
  The pairing $\blankeigpair{12}$ is nontrivial for $p = 691$ if and
  only if $\gr^{12} \mathfrak{h} \subsetneq \gr^{12} \g$.  In this
  case, there is a relation in $gr^{12} \g$,
  \begin{equation} \label{commutators}
    [\sigma_3,\sigma_9] \equiv 50[\sigma_5,\sigma_7]
    \pmod{691}.
  \end{equation}
\end{theorem}

\begin{proof}
  According to \cite{ihara:1999}, we have $\lambda_m(D_m) =
  1,2,16,144$ for $m = 3,5,7,9$ respectively.  If $\gr^{12}
  \mathfrak{h} \subsetneq \gr^{12} \g$, then the injectivity of
  $\iota$ implies that $\delta$ is contained in $\iota(\gr^{12} \g)$.
  By Theorem~\ref{notgenerated} and Lemma~\ref{Dmcondition}, equation
  \eqref{derivation12} becomes
  \begin{equation} \label{galois12}
    691c \cdot \bar{x}_{12} \equiv 190[\sigma_3,\sigma_9] +
        174[\sigma_5,\sigma_7]
    \pmod{691[\g,\g]}
  \end{equation}
  in $\gr^{12} \mathfrak{g}$ for some $c \not\equiv 0 \pmod{691}$.
  This yields \eqref{commutators}.  By the linear independence of
  $[\sigma_3,\sigma_9]$ and $[\sigma_5,\sigma_7]$, this equation must
  agree with that of \eqref{therelation} (noting
  Proposition \ref{filtrationlemma}, and after multiplication by an
  appropriate scalar).  Therefore, by Theorem~\ref{lierelation}, the
  pairing $\blankeigpair{12}$ is nontrivial.

  On the other hand, the coefficients in equation \eqref{galois12} of
  $[\sigma_i,\sigma_{12-i}]$ equal the values of the pairing
  $e_{i,12}$ up to a constant scalar multiple by the computation
  yielding Theorem~\ref{uniquepairing}.  Hence, by
  Theorem~\ref{notgenerated}, nontriviality of $\blankeigpair{12}$
  implies that $\gr^{12} \mathfrak{h} \subsetneq \gr^{12} \g$.
\end{proof}

In general, we expect that for $m$ minimal such that $p$ divides
$B_m$, a determination of the structure of $\gr^i \McD$ with $i
\le m$, together with a computation of a unique possibility for
the pairing $\blankeigpair{m}$ on $\mathcal{C} \times \mathcal{C}$
up to possibly trivial scalar (to show $e_{i,m} \neq 0$ for some
$i \in M_{\mathrm{o}}$ with $i < m/2$ if
$\eigpair{\mathcal{C}}{\mathcal{C}}{m} \neq 0$), would yield
(in essence) that Conjectures~\ref{k2conjecture},~\ref{iotaconj},
and the statement that $\gr^m \g \subsetneq \gr^m \h$, are
equivalent as well.

\section{Relationship with Greenberg's conjecture} \label{greenberg}

Let $K$ be a number field, $K_{\infty}$ the compositum of all
$\ZZ_p$-extensions of $K$, $\tilde{\Gamma} = \Gal(K_{\infty}/K)$
and $\tilde{\Lambda} = \zp[[\tilde{\Gamma}]]$.  Let $L_{\infty}$
be the maximal abelian unramified pro-$p$ extension of
$K_{\infty}$, and let $M_{\infty}$ be the maximal abelian
$p$-ramified pro-$p$ extension of $K_{\infty}$.  Then $X_{\infty}
= \Gal(L_{\infty}/K_{\infty})$ and $Y_{\infty} =
\Gal(M_{\infty}/K_{\infty})$ are $\tilde{\Lambda}$-modules.  We
say that a $\tilde{\Lambda}$-module is pseudo-null if its
annihilator has height at least $2$. Greenberg has made the
following conjecture \cite[Conjecture~3.5]{greenberg:2001}.

\begin{conjecture}[Greenberg]
  $X_{\infty}$ is pseudo-null as a $\tilde{\Lambda}$-module.
\end{conjecture}

For a prime ideal $\p$ of $K$ lying above $p$, let $r_{\p}$ be the
integer such that the decomposition group $\tilde{\Gamma}_{\p}
\subset \tilde{\Gamma}$ is isomorphic to $\zp^{r_{\p}}$.  In
certain cases, Greenberg's conjecture has an equivalent form in
terms of the torsion in $Y_\infty$. For example, we have the
following theorem, which is \cite[Corollary 14]{mccallum:2001}
(see also \cite{lannuzel-nguyenquangdo:2000}).
\begin{theorem}[McCallum]
  Assume that $r_{\p} \ge 2$ for all primes ${\p}$ above $p$ and that
  $\mu_{p^{\infty}} \subset K_{\infty}$.  Then Greenberg's conjecture
  holds if and only if $Y_{\infty}$ is $\tilde{\Lambda}$-torsion free.
\end{theorem}

Let $\McG$ be the Galois group of the maximal $p$-ramified pro-$p$
extension of $K$ and let $\tilde{G} \subset \McG$ be the Galois
group of the same extension over $K_{\infty}$.  Let $I(\McG)$
denote the augmentation ideal of $\ZZ_p[[\McG]]$.  Then the module
$Z = H_0(\tilde{G},I(\McG))$ has torsion subgroup isomorphic to
that of $Y_{\infty}$ (see \cite[Theorem 10]{mccallum:2001}).

Consider a free presentation of $\McG$ as in \eqref{presentation}
with minimal sets of generators and relations $X$ and $R$ as in
Section~\ref{intro}.  The presentation gives rise to an exact
sequence \cite{nguyenquangdo:1984, mccallum:2001}
\begin{equation} \label{Zexactseq}
  0 \to \tilde{\Lambda}^s \xrightarrow{\phi} \tilde{\Lambda}^g \to Z
\to 0,
\end{equation}
where $s = |R|$ and $g = |X|$.  Let $h$ denote the $\zp$-rank of
$\tilde{\Gamma}$ (so $h \ge g-s$).  Write $\tilde{\Lambda} =
\ZZ_p[[T]]$, where $T = (T_1,\ldots,T_h)$ and $1+T_i$ is the
restriction of a generator $x_i \in X$.  Furthermore, we choose
the remaining elements $x_i \in X$ with $h+1 \le i \le g$ such
that the image of $x_i$ in $\McG^{\ab}$ is torsion.  If
$$
f = (f_1,\ldots,f_g) \in \tilde{\Lambda}^g,$$ we also use the
notation
$$
f = \sum_{i=1}^g f_i dx_i.
$$

We briefly describe the map $\phi$ in \eqref{Zexactseq} (see
\cite{nguyenquangdo:1984} for details).  Fox \cite{fox:1953}
defines a derivative (extended to pro-$p$ groups in
\cite{nguyenquangdo:1984})
$$
D \colon \ZZ_p[[\McF]] \to I(\McF)
$$
satisfying $D(x) = x-1$ for $x \in X$ and a certain non-abelian
Leibniz condition that yields, for example,
$$
D [x,y] = (1-xyx^{-1})(x-1)+ (x-[x,y])(y-1)
$$
and
$$
D x^q = \left( \sum_{i=0}^{q-1} x^i \right) (x-1)
$$
for $q$ a power of $p$.  We also have the map
$$
\theta \colon I(\McF) = \sum \ZZ_p[[\McF]](x_i-1) \to \sum
\tilde{\Lambda} dx_i = \tilde{\Lambda}^g.
$$
We remark that the map $\Phi = \theta \circ D$ factors through
$\ZZ_p[[\McF/\McF'']]$ (where $N' = [N,N]$ for a group $N$).  The
map $\phi$ in the presentation \eqref{Zexactseq} is then obtained
on basis elements of $\tilde{\Lambda}^s$ by identifying them with
elements of $R$, considering these as elements of $\McF$ via the
presentation \eqref{intro}, and then composing with $\Phi$.
Roughly speaking, the map $\phi$ is the Jacobian matrix of the
relations.

We choose $R$ such that each $r_k \in R$ has the form
\begin{equation} \label{relation}
  r_k \equiv \left( \prod_{1\le i<j\le g} [x_i,x_j]^{f^k_{i,j}} \right) \cdot
  x_k^{l_k} \pmod{\McF''},
\end{equation}
with $g-s+1 \le k \le g$, where $l_k$ is a power of $p$ for $k >
h$ and $l_k = 0$ for $k \le h$ and where $f^k_{i,j} \in
\zp[[\McF^{\ab}]]$.  Each $f^k_{i,j}$ may be chosen not to involve
any $x_a$ with $i < a < j$ in order to make the expression unique.
Via the surjection $\zp[[\McF^{\ab}]] \to \tilde{\Lambda}$, we
obtain elements $f^k_{i,j}(T)$ of $\tilde{\Lambda}$, in which
$x_i$ is replaced by $1+T_i$ if $1 \le i \le h$ and $1$ otherwise.
We find
\begin{equation} \label{relimage}
  \Phi(r_k) = \sum_{1 \le i<j \le h} f^k_{i,j}(T)(-T_j dx_i+T_i dx_j)
  + \sum_{1 \le i \le h < j \le g} f^k_{i,j}(T) T_i dx_j + l_k dx_k.
\end{equation}

Let us assume from now on that $p$ is odd.  We require the
following lemma.
\begin{lemma} \label{UniversalNorm}
  Let $K$ be a number field containing $\mu_p$, and let $\alpha \in
  K^{\times}$ be a universal norm from the extension
  $K(\mu_{p^{\infty}})/K$.  Then the torsion in $\McG^{\ab}$ acts
  trivially on any $p$th root of $\alpha$.
\end{lemma}

\begin{proof}
  We may assume that $\alpha \notin \mu_{p^{\infty}}\cdot K^{\times
    p}$.  Let $K_n = K(\mu_{p^n})$, and let $\alpha_n \in
  K_n^{\times}$ be the elements of a norm-compatible sequence with
  $\alpha_1 = \alpha$.  Set
  \[
  \alpha'_n = \prod_{\sigma \in \Gal(K_n/K)}
  \sigma(\alpha_n)^{i_{\sigma,n}},
  \]
  with $i_{\sigma,n}$ the minimal nonnegative integer satisfying
  $\sigma(\zeta_{p^n}) = \zeta_{p^n}^{i_{\sigma,n}^{-1}}$.  Then
  $\alpha'_1 = \alpha$, and
  \[
  \alpha'_{n+1}\alpha'_n{}^{-1} \in K_{n+1}^{\times p^n}.
  \]
  Since $\alpha \notin \mu_{p^{\infty}} \cdot K^{\times p}$, we have
  that $\alpha'_n \notin K(\mu_{p^{\infty}})^{\times p}$.  By Kummer
  theory, the sequence $(\alpha'_n)$ defines a nontrivial element of
  $H^1(G, \displaystyle \lim_{\leftarrow} \mu_{p^n})$,
  with $G$ as in Section \ref{fundamental},
  on which
  $$
    \Gamma=\Gal(K(\mu_{p^{\infty}})/K)
  $$
  acts by the cyclotomic
  character.
  Hence, by a choice of isomorphism
  $$\lambda \colon \ZZ_p(1) \xrightarrow{\sim} \lim_{\leftarrow} \mu_{p^n}$$
  of $\McG$-modules and the fact that $G$ acts trivially on
  $\ZZ_p(1)$,
  it defines
  an element of $H^1(G,\ZZ_p)^{\Gamma}$.  Since $\Gamma$
  has cohomological dimension $1$, this element is the image under
  restriction of an element $\kappa \in \Hom(\McG,\ZZ_p)$.
  Furthermore, the map induced by the quotient $\ZZ_p \to \ZZ/p\ZZ$
  takes $\kappa$ to the Kummer character associated with $\alpha$,
  viewed as an element of $\Hom(\McG,\ZZ/p\ZZ)$ via $\lambda$ (as $K$
  contains $\mu_p$).  Since $\kappa$ has torsion-free image, the lemma
  is proven.
\end{proof}

We remark that all universal norms for
the cyclotomic
$\zp$-extension are
$p$-units.
We are now ready to prove the following consequence of
nontriviality of the pairing.

\begin{theorem} \label{greenbergconj}
  Let $K$ be a number field containing $\mu_p$ for an odd prime $p$,
  and let $S$ be the set
  of primes above $p$ and any real archimedean places. Suppose that
  $$\dim_{\ZZ/p\ZZ} H^2(G_{K,S}, \ZZ/p\ZZ) = 1.$$
  If the pairing
  $\Mcpairing{\alpha}{\beta}$ is nontrivial on two universal norms
  $\alpha, \beta$ for the extension $K(\mu_{p^{\infty}})/K$, then
  Greenberg's conjecture holds for $K$.
\end{theorem}

\begin{proof}
  By Lemma~\ref{UniversalNorm} and equation \eqref{BasicRelation}, the
  nontriviality of the cup product implies that
  $$
  f^g_{a,b}(0) \not\equiv 0 \pmod{p},
  $$
  with $f_{a,b}^{g}$ as in \eqref{relation} for some $a$ and $b$
  with $1 \le a < b \le h$.  Set $\mathfrak{m} =
  (p,T)\tilde{\Lambda}$.
  From equation \eqref{relimage}, we see that
  \begin{multline} \label{modrelimage}
    \Phi(r_g) \equiv \sum_{1 \le i<j \le h} f^g_{i,j}(0)(-T_j dx_i +
    T_i dx_j) \\ + \delta(l_g + \sum_{1 \le i \le g-1} f^g_{i,g}(0)T_i
    ) dx_g \pmod{\mathfrak{m}^2},
  \end{multline}
  where $\delta = 1$ if $h=g-1$ and $\delta = 0$ if $h=g$.  Equation
  \eqref{modrelimage} therefore shows that the coefficients of $dx_a$
  and $dx_b$ in $\Phi(r_g)$ are both nontrivial modulo
  $\mathfrak{m}^2$, with only the latter involving a linear term in
  $T_a$ that is nonzero modulo $p$.  Hence, one cannot factor a
  non-unit polynomial out of $\Phi(r_g)$, and $Z$ has no torsion.
\end{proof}

Let us specialize to the case $K = \QQ(\mu_p)$ by way of example.
Assume that $A_K$ has order $p$.  In \cite[Theorem
1]{mccallum:2001} (see also \cite{marshall:2001}), it is shown
that Greenberg's conjecture holds if $f(0) \not\equiv -pf'(0)
\pmod{p^2}$, where $f$ denotes a characteristic power series of
the $p$-part of the class group of $\QQ(\mu_{p^{\infty}})$.
Essentially, this is a condition on the last term in the
expression \eqref{modrelimage}.  (By the discussion in
Section~\ref{GaloisGroup}, it is exactly that $l_g \not\equiv 0
\pmod{p^2}$.)  Since that term involves only $dx_g$, the method of
analysis in Theorem~\ref{greenbergconj} seems incapable of giving
an alternate proof of McCallum's result. On the other hand, since
cyclotomic $p$-units are universal norms, we have the following
variation on McCallum's result which replaces the condition on the
characteristic power series with a condition on the pairing.

\begin{corollary} \label{greenbergcycl}
  Let $K = \QQ(\mu_p)$, and assume that
  $A_K$ is a cyclic group.
  If the restriction of $\blankpairing$ to the cyclotomic $p$-units is
  nontrivial, then Greenberg's conjecture holds for $K$.
\end{corollary}

Using Theorem~\ref{NontrivialPairing}, this gives another proof of
Greenberg's conjecture for $\mathbb{Q}(\mu_{37})$.

\providecommand{\href}[2]{#2}

\end{document}